\documentclass[11pt, leqno]{article}
\usepackage{times, amsfonts, amsmath, amssymb, latexsym, setspace, epsfig}
\usepackage{subfigure, graphics}

\def\bs{\boldsymbol}
\def\Ex{{\rm I\!E}}
\def\Pr{{\rm I\!P}}
\def\be{\begin{equation}}
\def\ee{\end{equation}}
\def\bea{\begin{eqnarray*}}
\def\eea{\end{eqnarray*}}
\def\bean{\begin{eqnarray}}
\def\eean{\end{eqnarray}}
\def\nn{\nonumber}
\def\nin{\noindent}

\def\ra{\rightarrow}

\def\Bl{\Bigl}
\def\Br{\Bigr}
\def\wt{\widetilde}

\def\R{{\bf R}}

\def\alp{\alpha}

\def\eps{\epsilon}

\def\th{\theta}

\def\Yn{\mathbf{Y}\!_n}
\def\Zn{\mathbf{Z}_n}

\def\logLRI{{\rm logLR}_I}
\def\logLRIIc{{\rm logLR}_{I,I^c}(\hat{\th}_{[n]})}

\newtheorem{Theorem}{Theorem}
\newtheorem{Proposition}{Proposition}

\newtheorem{Lemma}{Lemma}

\begin{document}
\setstretch{1.1}

\title{Calibrating the scan statistic: finite sample performance vs. asymptotics}

\author{Guenther Walther\thanks{Work supported by NSF grants DMS-1501767 and DMS-1916074} 
        and Andrew Perry\thanks{Work supported a Summer Undergraduate Research Grant from the
Vice Provost for Undergraduate Education at Stanford University}\\
        Department of Statistics, 390 Serra Mall \\
        Stanford University, Stanford, CA 94305 \\
        gwalther@stanford.edu}
\date{July 2021}
\maketitle

\begin{abstract}
We consider the problem of detecting an elevated mean on an interval with unknown location
and length in the univariate Gaussian sequence model. 
Recent results have shown that using scale-dependent critical values for the scan statistic
allows to attain asymptotically optimal detection simultaneously for all signal lengths, 
thereby improving on the traditional scan, but this procedure has been criticized
for losing too much power for short signals. We explain this discrepancy by showing  that these
asymptotic optimality results will necessarily be too imprecise to discern the performance
of scan statistics in a practically relevant way, even in a large sample context. 
Instead, we propose to assess the performance 
with a new finite sample criterion. We then present three calibrations for scan statistics
that perform well across a range of relevant signal lengths:
The first calibration uses a particular adjustment to the critical values and is therefore
tailored to the Gaussian case. The second calibration uses a scale-dependent adjustment to the 
significance levels and is therefore applicable to arbitrary known null distributions. The third
calibration restricts the scan to a particular sparse subset of the scan windows and then
applies a weighted Bonferroni adjustment to the corresponding test statistics. This 
{\sl Bonferroni scan} is also applicable to arbitrary null distributions and in addition is very simple to 
implement. We show how to apply these calibrations for scanning in a number of distributional settings:
for normal observations with an unknown baseline and a known or unknown constant variance,
for observations from a natural exponential family,
for potentially heteroscadastic observations from a symmetric density by employing self-normalization
in a novel way, and for exchangeable observations using tests based on permutations, ranks or signs.
\end{abstract}

\vfill
\vfill



\newpage

\section{Introduction}  \label{introduction}

There has been an considerable amount of recent work
that requires a good solution to the following problem: We observe
\be  \label{model}
Y_i\ =\ f_n(i) + Z_i,\ \ \ \ i=1,\ldots,n,
\ee
where the $Z_i$ are i.i.d. $N(0,1)$ and $f_n(i)=\mu_n\,{\bs 1}(i \in I_n)$ with 
$I_n=(j_n,k_n]$, $0\leq j_n<k_n\leq n$. The task is to decide whether a signal is present,
i.e. to test $\mu_n=0$ vs. $\mu_n >0$, when both the amplitude $\mu_n$ and the support $I_n$
are unknown. This is the prototypical model for the problem of detecting objects in noisy data
and for certain goodness-of-fit tests, see the references below. The standard approach to this
problem is based on the {\sl scan statistic} (generalized likelihood ratio statistic)
\be \label{scan}
Scan_n(\Yn)\ =\ \max_I T_I(\Yn),
\ee
where for intervals $I=(j,k]$, $0\leq j<k\leq n$, we write
\be \label{T}
T_I(\Yn):= \frac{\sum_{i \in I} Y_i}{\sqrt{|I|}} = \frac{\sum_{i =j+1}^k Y_i}{\sqrt{k-j}},
\ee

see e.g. Glaz et al. (2001) or Arias-Castro et al. (2005). Naus and Wallenstein~(2004) observed
that the scan statistic is sensitive for the detection of signals with very short support $|I_n|$
at the expense of signals with moderate and large support. A heuristic explanation of this effect
is as follows: There are about $\frac{n}{w}$ disjoint intervals $I$ with length $|I|=w$. Since the
corresponding $T_I(\Zn)$ are independent $N(0,1)$, their maximum concentrates around 
$\sqrt{2 \log \frac{n}{w}}$. This remains true for the maximum over all (overlapping)
intervals of length $w$
because the corresponding $T_I(\Zn)$ are strongly correlated with  the set of $\frac{n}{w}$
independent ones. Hence the distribution of $Scan_n(\Zn)$ is dominated by the small intervals
with $|I| \approx 1$ and concentrates around $\sqrt{2 \log n}$; see Siegmund and Venkatraman (1995)
for a formal proof. This heuristic suggests that the scan does not provide an optimal aggregation
of the evidence for the various interval lengths $w$, and that there may be a `free lunch': since
$\max_{I:|I|=w} T_I(\Zn)$ concentrates around $\sqrt{2 \log \frac{n}{w}}$, it may be possible
to increase $T_I(\Yn)$ by $\sqrt{2 \log n}-\sqrt{2 \log \frac{n}{|I|}}$ without noticeably
changing the null distribution of $\max_I T_I(\Zn)$, thereby increasing the power at larger scales
and remedying the problem described by Naus and Wallenstein~(2004). This idea was formalized
by D\"{u}mbgen and Spokoiny~(2001) who introduced the statistic
\be  \label{DS}
DS_n(\Yn)\ =\ \max_I \left(T_I(\Yn) -\sqrt{2\log \frac{en}{|I|}}\right)
\ee
and established various optimality results. These results show
that $DS_n$ leads to asymptotically optimal detection in the model (\ref{model}) for all
scales $|I_n|$, in the sense
that no other statistical test can improve on $DS_n$ even if
the scale $|I_n|$ of the signal were known in advance. In other words, while scanning over locations
leads to an unavoidable multiple testing penalty,
there is no further material price to pay in the asymptotic minimax framework
for scanning over multiple scales $|I|$ when using $DS_n$. In contrast, it was shown in
Chan and Walther~(2013) that inference using $Scan_n$ will be suboptimal on all but the smallest scales.

The asymptotic optimality of $DS_n$ across all scales has made this statistic a popular choice
for a range of problems involving the detection of signals or testing goodness-of-fit, see e.g. 
D\"{umbgen} and Walther~(2008), Rohde~(2008), Frick et al.~(2014), K\"{o}nig et al.~(2020), or
Datta and Sen~(2018).
However, $DS_n$ has been criticized by Siegmund~(2017) for losing too much detection power
on small scales. Indeed, Figure~\ref{fig1} shows
 a plot of the critical values for $T_I(\Yn)$ as a function
of $|I|$, for various methods for aggregating these statistics when $n=10^4$. 
The green line shows\footnote{The procedures in Sections \ref{methods} and \ref{sims} 
use interval lengths $|I|$ that are bounded by $\frac{n}{4}$, but the critical values are similar
and the conclusions described here are not 
sensitive to that upper bound.}
 the critical values
resulting from the calibration $DS_n$. Compared to the black line resulting from $Scan_n$, it is clear
that the price for obtaining smaller critical values (and hence more power) at larger scales comes at the
price of larger critical values (and thus reduced power) at small scales. For example, if $n=10^3$
then the exact critical
value at $|I|=1$ is 4.14 with $Scan_n$ and 5.09 with $DS_n$. So in order to declare a discovery,
$T_I(\Yn)$ needs to clear the 4$\sigma$ threshold when using the calibration $Scan_n$, but it
needs to exceed the 5$\sigma$ threshold with the calibration $DS_n$. This is arguably an unacceptable
a price to pay if the concern is about the discovery of a signal on a small scale $|I_n|$,
which is typically the case in applications using the model (\ref{model}), as pointed out
by Siegmund~(2017). It is also informative to examine this discrepancy from the angle of the
multiple testing problem which lies at the heart of model (\ref{model}).
Looking at different sample sizes, one finds that the critical
value for $Scan_n$ is 4.71 when $n=10^4$ and 5.21 for $n=10^5$.
Thus the penalty incurred by using $DS_n$ on the small scales is equivalent to the multiple testing
penalty incurred when increasing the sample size by a factor of between 10 and 100.
Furthermore, simulations and theoretical considerations show that this discrepancy between $DS_n$ and
$Scan_n$ will not disappear asymptotically. Therefore, these finite sample results 
clearly favor $Scan_n$ over $DS_n$, despite the strong theoretical support for the latter.

In Section~\ref{optimality} we explain why the practical performance of $DS_n$ can be markedly 
inferior to that of $Scan_n$ despite the asymptotic optimality results for $DS_n$. 
We introduce an exact finite sample criterion that is a more useful measure of the performance of
a scan statistic. We then
present three calibrations that also satisfy these optimality results (and hence allow detection
thresholds that are substantially lower than those for $Scan_n$ at larger scales) without
paying a material price at smaller scales. The first calibration uses a particular adjustment to the
critical values of the $T_I(\Yn)$ and is therefore tailored to the Gaussian case in (\ref{model}).
The second calibration involves an adjustment to the significance levels of the $T_I(\Yn)$
and is therefore not specific to the Gaussian case but can be applied to arbitrary null distributions
of $\Zn$. The third calibration restricts the scan to a particular sparse approximation set of
intervals $I$ and then simply applies a weighted Bonferroni adjustment to those $T_I(\Yn)$. This
calibration is also applicable for arbitrary null distributions of $\Zn$ and furthermore has the
advantage that the adjustment via Bonferroni removes the need to approximate critical values
by simulation or by analytical approximation. For each of these three calibrations, the key
to achieving the desired performance is a slight relenting on the demand for optimality
 at the largest scales $|I_n| \approx n$,
which is arguably not of practical concern. Section~\ref{sims} examines the theoretical and
practical performance of these calibrations in terms of the new finite sample criterion. 
In Section~\ref{epidemic} we show how to apply these calibrations for scanning in a number
of distributional settings:  
for normal observations with an unknown baseline and a known or unknown constant variance,
for observations from a natural exponential family,
for potentially heteroscadastic observations from a symmetric density by employing self-normalization
in a novel way,
and for exchangeable observations using tests based on permutations, ranks or signs.
Section~\ref{other} briefly discusses
how these calibrations extend to other settings such as observations from densities and
the multivariate case. Proofs are deferred to the Appendix.

\section{Asymptotic optimality and finite sample behavior}  \label{optimality}

Siegmund and Venkatraman~(1995) show that
\be  \label{SV}
Scan_n(\Zn) \ =\ \sqrt{2 \log n} + O_p\Bl(\frac{\log \log n}{\sqrt{\log n}}\Br)\ =\ 
\sqrt{2 \log n} +o_p(1)
\ee
A heuristic explanation of this result is as follows: While $Scan_n(\Zn)$ is the maximum of $\sim n^2$
i.i.d. standard normals, the first order term $\sqrt{2 \log n}$ shows that it behaves likes the
maximum of only $\sim n$ i.i.d. standard normals, 
the reason being that many $T_I(\Zn)$ are correlated since the $I$ are overlapping.
Arias-Castro et al. (2005) show that  if
$\sqrt{|I_n|} \,\mu_n \geq \sqrt{(2+\eps)\log n}$ for
$\eps >0$, then $Scan_n$ will be {\sl asymptotically powerful}, i.e. the probabilities of type I
and of type II error both go to zero asymptotically. They generalize the setting (\ref{model})
to multivariate and geometrically defined signals and show that these detection problems
give rise to similar detection boundaries of the form $\sqrt{2D \log n}$, which determine
the amplitude of the signal that is detectable. 
To appreciate the importance of the constant $D$, note that $\sqrt{2 \log n}$ is essentially
the Bonferroni adjusted critical value for $n$ independent z-tests.
Hence the threshold $\sqrt{2D \log n} =\sqrt{2 \log n^D}$ corresponds
to a multiple testing problem whose difficulty is determined by $n^D$ independent z-tests.
For this reason Arias-Castro et al. (2005) call $D$ the exponent of effective dimension.

A key insight of 
D\"{u}mbgen and Spokoiny~(2001) is that employing critical values that depend on the scale 
$|I_n|$ allows to detect even smaller amplitudes, so the above detection boundary can in
fact be improved upon: Translating their methodology into the setting (\ref{model}) one can show 
that if 
\be  \label{consist}
\sqrt{|I_n|} \,\mu_n\ \geq \ \sqrt{(2+\eps_n)\log \frac{n}{|I_n|}}
\ee
with $\eps_n \ra 0$ not too fast, namely $\eps_n \sqrt{\log \frac{n}{|I_n|}} \ra \infty$,
then $DS_n$ has asymptotic power 1. Thus, if e.g. $|I_n|/n=n^{-1/2}$, 
then the critical multiplier for $\log n$ in the detection boundary can be reduced from 
2 to 1 if one
uses $DS_n$ in place of $Scan_n$. This multiplier  cannot be reduced further,
as reliable detection becomes asymptotically impossible for any procedure if `$2+\eps_n$'
is replaced by `$2-\eps_n$' in (\ref{consist}), 
as can be seen from the lower bounds given by D\"{u}mbgen and Spokoiny~(2001)
in the case of small scales $|I_n|$ and by D\"{u}mbgen and Walther~(2008) in the case of large
scales.

 While the above results provide a precise characterization of the detection boundary,
we will now argue that these asymptotic results will necessarily be too imprecise to discern
practically relevant performance characteristics
 of these max-type statistics with an appropriate level of precision, even in the
large sample context. This is illustrated with the following theorem:

\begin{Theorem}   \label{thm1}
If there exists a sequence of critical values $\{\kappa_n \}$ for which $Scan_n$ is 
asymptotically powerful
against $\{f_n\}$ in the model (\ref{model}) with $|I_n| \leq n^p$ for $p<\frac{1}{16}$,
then $Scan_n$ is also asymptotically powerful for any sequence of critical values 
$\tilde{\kappa}_n = \kappa_n +O(1)$ with $O(1) \geq 0$.
\end{Theorem}

Hence if the critical values $\kappa_n$ result in an asymptotically powerful test, then so will
the critical values $\kappa_n +100$ (say), even though the latter are clearly not a useful
choice even if the sample size is enormous.
As explained in the proof of the theorem, analogous conclusions hold for related
test statistics such as $DS_n$, or when considering asymptotic optimality 
with a fixed significance level. The upshot of the theorem is that asymptotic
optimality statements such as (\ref{consist}) will characterize optimal critical 
values only to a precision
of $O(1)$, while the null distribution (\ref{SV}) concentrates with a rate of $o(1)$
around $\sqrt{2 \log n}$. Importantly, it is the $o(1)$ term that determines
relevant performance characteristics
of the statistic. To see this heuristically, note that the $\sqrt{2 \log n}$ term 
arises as a multiple testing adjustment for $\sim n$ independent test statistics. 
Suppose we greatly increase the multiple testing problem by a factor of $10^k$. Since
$$
\sqrt{2 \log (10^k n)} \leq \sqrt{2 \log n} + \frac{2k}{\sqrt{\log n}} = \sqrt{2 \log n}+o(1),
$$
one sees that the adjustment for this much larger multiple testing problem will affect the null
distribution only on the scale of the $o(1)$ term and will therefore be overlooked by the
asymptotic theory. But it is well known in statistical practice 
that a Bonferroni-like correction by such a large
factor will typically affect the power of the statistic in a way that is quite relevant for 
inference. 

In summary, while the current state of the asymptotic optimality theory allows to derive the fundamental 
difficulty of the detection problem in terms of the optimal detection threshold, such as (\ref{consist}),
and this provides a necessary condition for the large sample optimality of a test statistic,
the above considerations show that these conditions are not sufficient for a good performance
of the test statistic, even in the large sample setting.
This raises the question of how one should evaluate the performance of calibrations such as $Scan_n$ 
and $DS_n$ in a practically informative way. One option would be to develop a refined asymptotic
theory. But there are reasons to doubt whether this would be informative. For example,
D\"{u}mbgen and Spokoiny~(2001) introduce a refined version of the statistic $DS_n$ which
employs an iterated logarithm, see also Proksch et al. (2018). 
An inspection of this refinement suggests that it should improve
the performance on small scales, but our simulations show that it is nearly indistinguishable
from $DS_n$ for sample sizes up to $10^6$, presumably because the asymptotics set in too
slowly. For this reason we suspect that a refined asymptotic optimality theory may likewise
not sufficiently illuminate the performance of a calibration. Therefore
we propose a different approach  which focuses on the finite sample performance
of the statistic: 

The idea is to define the {\sl realized exponent} $e_n(|I_n|)$ as the solution of the equation
\be \label{realizedexp}
\sqrt{|I_n|}\,\mu_{min}(n,|I_n|)\ =\ \sqrt{2 e_n(|I_n|) \log \frac{e n}{|I_n|}}
\ee
where $\mu_{min}(n,|I_n|)$ as the smallest change in mean that a calibration
is able to detect reliably. In more detail, in model (\ref{model}) we are looking for the
the smallest $\mu_n$ on a randomly placed interval $I_n$ with given length $|I_n|$ that the
test statistic will detect with power 80\% at the 10\% significance level.
Note the formal similarity of $e_n(|I_n|)$ to the exponent of effective dimension $D$ described above, which
measures the fundamental difficulty of the detection problem.
But in contrast to the exponent of effective dimension, the realized exponent $e_n(|I_n|)$
is a function of the test statistic and shows how
close the test statistic comes to attaining the asymptotic detection boundary 
in the finite sample situation at hand. Importantly, $e_n(|I_n|)$ can be readily computed
via Monte Carlo, see Section~\ref{sims}. $e_n(|I_n|)$ is therefore well suited to
compare various test statistics to the benchmark given by the traditional and popular
statistic $Scan_n$. 

The asymptotic detection threshold (\ref{consist}) suggests that
a statistic that performs well for the detection problem (\ref{model}) should have
a realized exponent $e_n(|I_n|)$ close to 1 for all scales $|I_n|$. As an example,
for $n=10^6$ we find that $Scan_n$ has a realized exponent $e_n(1)=1.32$ while 
$e_n(n^{\frac{1}{2}})=2.45$, see Section~\ref{sims}.
The poor performance at the scale $|I_n|=n^{\frac{1}{2}}$ is not surprising since
it was shown in Chan and Walther~(2013) that $Scan_n$ will (asymptotically) attain the optimal 
threshold (\ref{consist}) only at the smallest scales. $DS_n$ is designed to attain
asymptotic optimality across all scales. Its realized exponent
$e_n(n^{\frac{1}{2}})=1.87$ shows indeed a considerable improvement over $Scan_n$, but the
price for this improvement is the disappointing performance at the important small scales:
 $e_n(1)=1.67$. This means that the power loss of $DS_n$ as compared to
$Scan_n$ is equivalent to increasing the size of the multiple testing problem by a factor of 
$n^{1.67-1.32} \approx 126$, validating the criticism of $DS_n$ that was stated in the Introduction.

This immediately raises the question whether this outcome represents an unavoidable
trade-off, or whether it is possible to construct a statistic
that attains the advantageous performance of $DS_n$ at larger scales without sacrificing
performance at small scales. In the next section we will answer this question to the affirmative
by presenting three different approaches for calibrating the various scales of such 
max-type statistics: the first approach calibrates the critical values,
the second calibrates the significance levels, and third one introduces a weighted
Bonferroni scheme on a sparse subset of the scan windows which is very simple to implement.

\section{Three ways of calibrating scan statistics for good performance}  \label{methods}

For each of the following three calibrations, the key to achieving a good finite sample
performance is to give up some power at the largest scales $|I_n| \approx n$. As explained
in the Introduction, those scales are typically not of a concern, and it turns out that 
a rather small sacrifice in power there will produce a considerable improvement in finite sample 
performance.

\subsection{A variation of the Sharpnack-Arias-Castro calibration}

\begin{figure}
\centering
\subfigure{\includegraphics[width=0.8\textwidth]{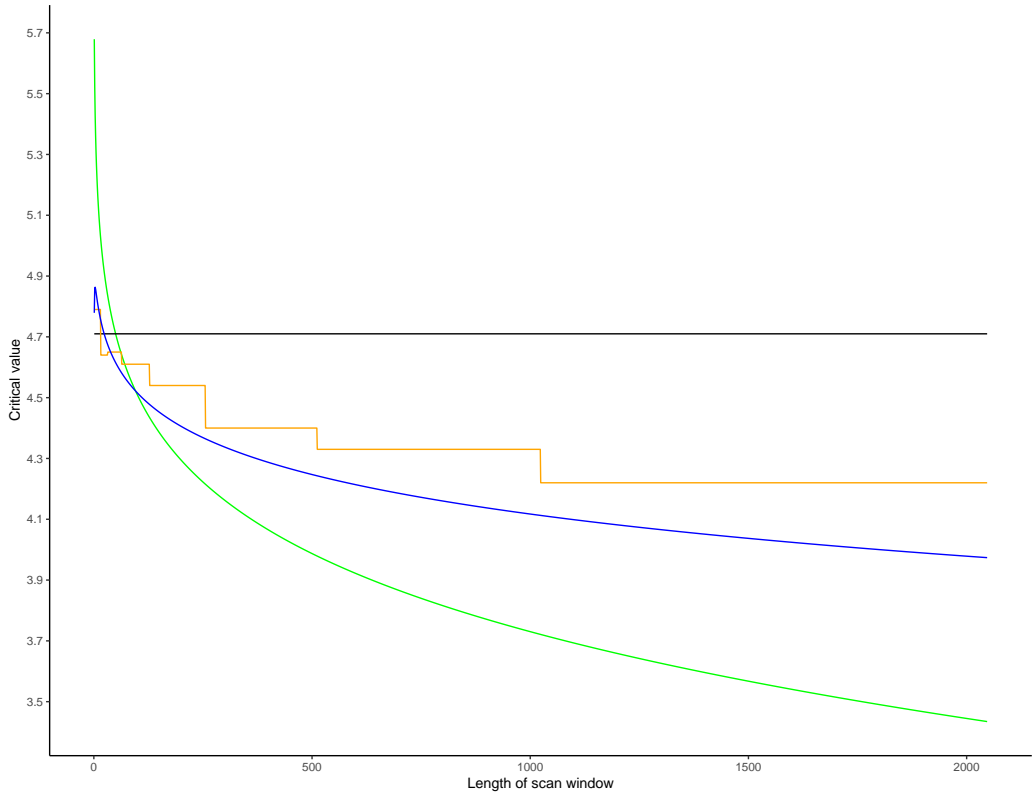}} 
    \subfigure{\includegraphics[width=0.8\textwidth]{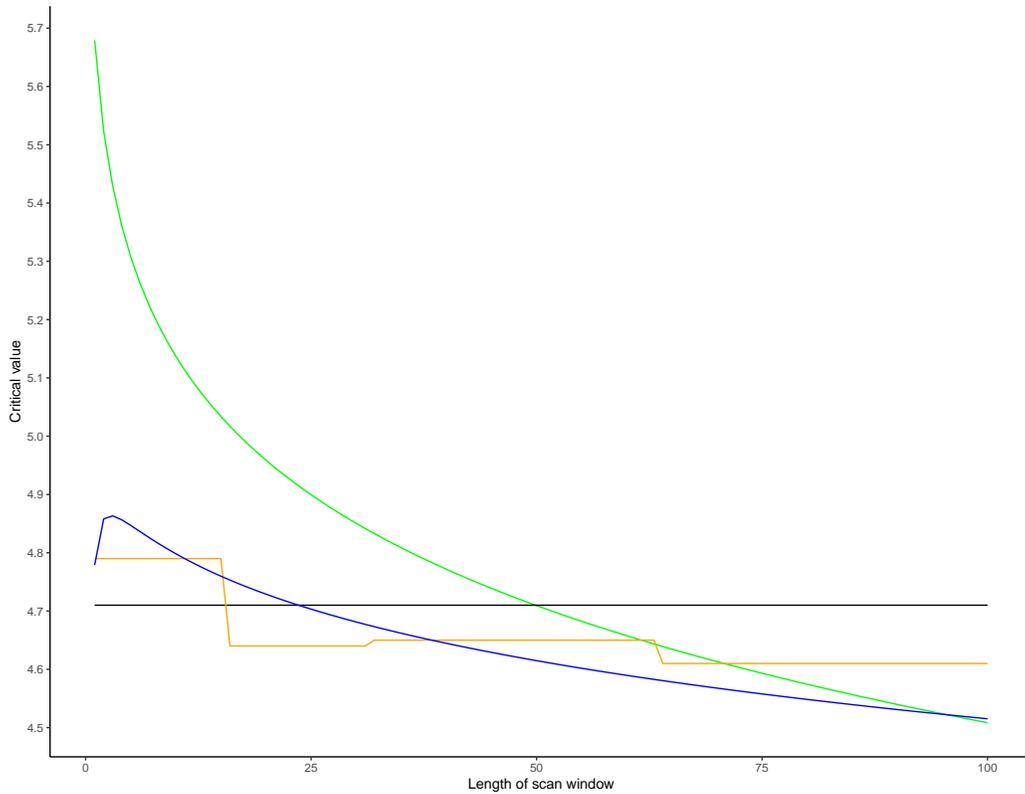}}
\vspace{-0.5cm}
\caption{Critical values as a function of the length of the scan window
for $Scan_n$ (black line), $DS_n$ (green), $SAC_n$ (blue) and
{\sl BlockedScan$_n$} (orange).
Sample size is $n=10^4$ and significance level is 10\%. The critical values were simulated
with $10^4$ Monte Carlo simulations, using the same largest window length of about
$\frac{n}{4}$ for all four procedures. The bottom plot zooms in on the smallest window lengths.}
\label{fig1}
\end{figure}

The first calibration is a simplified version of a standardization used by Sharpnack
and Arias-Castro~(2016) in the context of proving a limiting distribution. We therefore
call this correction to the critical values of the scan the {\sl Sharpnack-Arias-Castro calibration}
\be  \label{SAC}
SAC_n(\Yn)\ =\ \max_I \left(T_I(\Yn) -\sqrt{2\log \Bl[\frac{en}{|I|}(1+\log |I|)^2\Br]}\right).
\ee
This statistic is similar to $DS_n$ given in (\ref{DS}), but the factor $(1+\log |I|)^2$
increases the penalty for larger scales and therefore transfers some power to smaller scales.
The resulting critical values $T_I(\Yn)$ are therefore
$$
\sqrt{2\log \Bl[\frac{en}{|I|}(1+\log |I|)^2\Br]} + q_n(\alpha) 
$$
where $q_n(\alpha)$
is the $(1-\alp)$-quantile of $SAC_n(\Zn)$, which is obtained by simulation, see Section~\ref{sims}.
These critical values are plotted as a function of $|I|$  in Figure~\ref{fig1} (blue line) for sample 
size $n=10^4$. For comparison, the plot also shows the critical values resulting from the
calibration $DS_n$ (green line), which uses the correction term $\sqrt{2\log \frac{en}{|I|}}$
added to the $(1-\alp)$-quantile of $DS_n(\Zn)$, as well as those resulting from $Scan_n$ (black line),
which do not depend on the scale $|I|$ and are simply given by the $(1-\alp)$-quantile of 
$Scan_n(\Zn)$. The plot shows that while $DS_n$ has critical values that are
substantially smaller than those of $Scan_n$ for most of the scales $|I|$ (and therefore
$DS_n$ has more power there), it also has considerably larger critical values and
hence inferior performance at the smallest scales. In contrast, the critical values of
$SAC_n$ are only slightly larger than those of $Scan_n$
at the smallest scales while still being considerably
smaller at the other scales. Therefore the simple correction term in $SAC_n$
arguably produces a useful improvement over the traditional scan.



We note that the correction term in $SAC_n$ is specific to the Gaussian tails of the test statistic
$T_I(\Zn)$. In order to use $SAC_n$ for different settings it is therefore necessary to
transform the statistic such that it has (sub)Gaussian tails under the null hypothesis.
While this is possible in many cases, e.g. Rivera and Walther~(2013) and Walther~(2021) show how to do this
for likelihood ratio statistics by employing the square-root of the log likelihood ratio,
it is desirable to devise a calibration that is applicable to more general null distributions
of the test statistic. 

\subsection{The blocked scan}

To this end we propose to calibrate the significance level rather
than the critical value, following the idea of the {\sl blocked scan} introduced in
Walther~(2010). The blocked scan groups intervals of roughly the same length (e.g. having length
within a factor of 2) into a block. Then all intervals within a block are assigned the
same critical value. The significance level for each block is set so that the resulting
test performs well across all scales. It turns out that this can be achieved by
assigning a significance level that is proportional to a harmonic sequence. In more detail:

We define the first block to comprise the smallest intervals with length up to about $ \log n$,
namely $|I|<2^{s_n}$, where $s_n:=\lceil \log_2 \log n \rceil$. From there on we use powers of 2 to 
group interval lengths: The $B$th block contains all intervals with length $|I| \in
[2^{B-2+s_n},2^{B-1+s_n})$, $B=2,\ldots,B_{max}:=\lfloor \log_2 \frac{n}{4} \rfloor -s_n+1$.
The choice of $B_{max}$ means that the largest intervals we consider have length about $\frac{n}{4}$;
we found that the results in this paper are not sensitive to this endpoint and all simulations
reported in this paper use this endpoint for the maximum in (\ref{scan}),(\ref{DS}),(\ref{SAC}).
Now we let the significance level of the $B$th block decrease like a harmonic sequence:
$BlockedScan(\Yn)$ rejects if for any block $B$:
$$
\max_{I \in B\mbox{th block}} T_I(\Yn) > c_{B,n} \left(\frac{\tilde{\alp}}{B}\right)
$$
where $c_{B,n}(\alp)$ is the $(1-\alp)$-quantile of $\max_{I \in B\mbox{th block}} T_I(\Zn)$,
which is obtained by simulation, and $\tilde{\alp}$ is set so that 
$BlockedScan(\Zn)$ has overall level $\alpha$, see Walther~(2010). 

Letting the significance level of the $B$th block {\sl decrease} with $B$ like a harmonic sequence,
rather than letting it {\sl increase} as in the original proposal for the blocked scan in Walther~(2010),
has the same effect as using the correction term $SAC_n$ instead of $DS_n$ for the critical values:
Asymptotic optimality at the largest scales $|I| \approx n$ will be lost, but one gains a notably better
performance at the smallest scales. This can be seen in Figure~\ref{fig1}, which shows that the
critical values $c_{B,n}$ (orange line) mimic those of $SAC_n$.

\subsection{The Bonferroni scan}

Note that the critical values of all four statistics considered so far need to be approximated
either analytically or by simulation. Our third proposal avoids the effort that comes with
 such an approximation if the
null distribution of $T_I(\Zn)$ is known, as is the case for model (\ref{model}) where it 
is standard normal. In order to explain the idea it is helpful to first review how the critical
values of the above four statistics can be approximated.
 Note that the scan in (\ref{scan}) is defined as the 
maximum over $\sim n^2$ intervals. Hence the simulation of the null distribution becomes
computationally infeasible for larger $n$, which has motivated analytical approximations to
quantiles such as in Siegmund and Venkatraman (1995). Alternatively, work
developed in the last 15 years
has shown that one can effectively approximate this maximum by evaluating
$T_I$ only on a judiciously chosen set of intervals $I$. Importantly, it suffices to use
only about $O(n)$ intervals, hence critical values can be readily simulated via Monte Carlo. 
Our third proposal exploits the sparsity of such a collection of intervals not only for computation,
but also for inference. The idea is that if the collection is sparse enough, then a simple
union bound may produce critical values for the local test statistics that will result in a good
performance. Note that this requires to strike a delicate balance in order to guarantee good power:
the collection of intervals has to be rich
enough so it can provide a good approximation for an arbitrary interval $I$,
but it also has to be sparse enough so that a Bonferroni adjustment will not unduly diminish the power 
of the local statistics. We will demonstrate below that applying a weighted Bonferroni adjustment
to the approximating set of intervals introduced 
by Walther~(2010) (see Rivera and Walther~(2013) for the univariate version used here) results in
a test that does
indeed perform nearly as well as $SAC_n$ and $BlockedScan_n$. We call this calibration
the {\sl Bonferroni scan}.

For integers $\ell \geq 0$ define $m_{\ell}:=2^{\ell}$ and $d_{\ell}:=\left\lceil m_{\ell}/
\sqrt{2\log \frac{en}{m_{\ell}}} \right\rceil$. Then we approximate intervals with lengths in
$(m_{\ell},2m_{\ell}]$ with intervals from the collection
$$
{\cal J}_{\ell}:= \Bl\{ (j,k]: j,k \in \{i d_{\ell}, i=0,1,2,\ldots\} \mbox{ and }
m_{\ell} \leq k-j < 2m_{\ell} \Br\}
$$
which is essentially the collection given in Rivera and Walther~(2013) but indexed backwards.
The spacing $d_{\ell}$ implies that the approximation error relative to the length of the interval
becomes smaller for smaller intervals at the rate $(2 \log (en/m_{\ell}))^{-1/2}$, and this 
rate guarantees
both a good enough approximation as well as a sufficiently sparse representation.
This is an important difference to other approximation schemes introduced in the literature,
such as the ones given in Arias-Castro et al. (2005)
or in Sharpnack and Arias-Castro~(2016), although it may be possible to modify those to
produce comparable results.

Now we proceed as with the blocked scan above: We group intervals into blocks and then assign
a significance level that is proportional to a harmonic sequence. So we define
$$
\mbox{$B$th block }:= \left\{ \begin{array}{ll}
  \bigcup_{\ell =0}^{s_n-1} {\cal J}_{\ell} & \mbox{ if $B=1$}\\
  {\cal J}_{B-2+s_n} & \mbox{ if $B=2,\ldots,B_{max}$}
  \end{array}
\right.
$$
Hence the $B$th block contains intervals of exactly the same lengths as for the blocked
scan, but the endpoints of the intervals are thinned out  with the spacing $d_{\ell}$.
In contrast to the blocked scan, we can directly find the critical value for the $B$th
block by using a Bonferroni adjustment: If $v_n(\alp)$ denotes the $(1-\alp)$ quantile
of $T_I(\Zn)$, then the critical value for the $B$th block is
$$
v_n\left(\frac{\alp}{\# (B\mbox{th block})\, B \,\sum_{i=1}^{B_{max}} \frac{1}{i}} \right)
$$
and {\sl BonferroniScan$_n$}$(\Yn)$ rejects if for any block $B$ the statistic
$\max_{I \in B\mbox{th block}} T_I(\Yn)$ exceeds that critical value. 
The number of intervals in the $B$th block, $\# (B\mbox{th block})$, provides the
Bonferroni adjustment within each block, while the weighted Bonferroni
adjustment across blocks is given by $B^{-1}$, standardized by $\sum_{i=1}^{B_{max}} \frac{1}{i}$.
Therefore it follows immediately that this test has level at most $\alp$.

\begin{figure}
\centering
\subfigure{\includegraphics[width=0.8\textwidth]{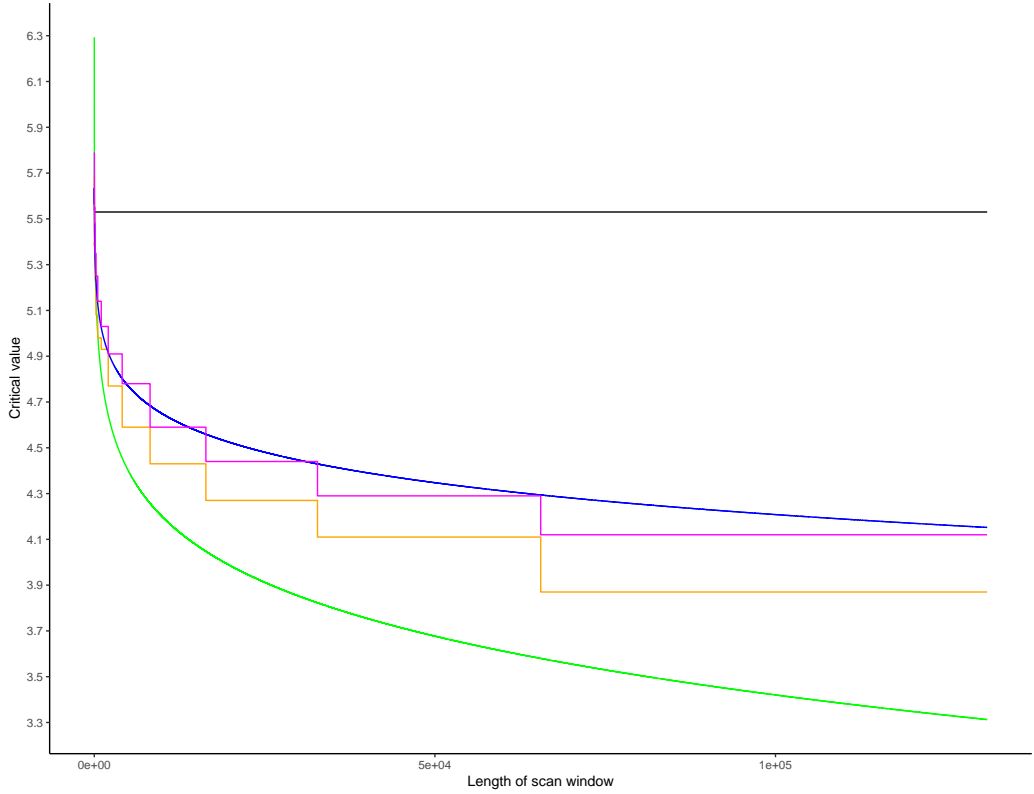}}
    \subfigure{\includegraphics[width=0.8\textwidth]{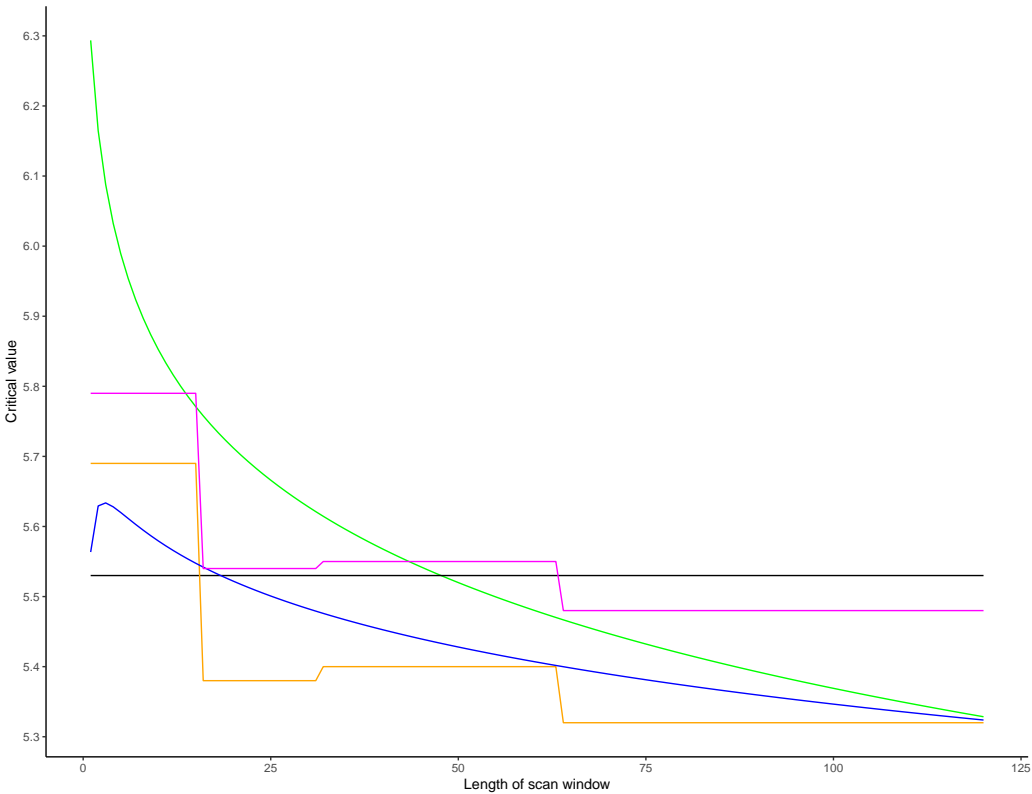}}
\vspace{-0.5cm}
\caption{Critical values as a function of the length of the scan window
for $Scan_n^{app}$ (black line), $DS_n^{app}$ (green), $SAC_n^{app}$ (blue), 
{\sl BlockedScan$_n^{app}$} (orange) and the Bonferroni scan (magenta). 
Sample size is $n=10^6$ and significance level is 10\%. The critical values were simulated
with $10^4$ Monte Carlo simulations, using the same largest window length of about
$\frac{n}{4}$ for all four procedures. The bottom plot zooms in on the smallest window lengths.}
\label{fig2}
\end{figure}

Of course one may also evaluate all of the above calibrations on this approximating set of intervals.
We denote the resulting procedures with the superscript $^{app}$. Using this approximating set 
allows for a much faster evaluation of the statistic as well as simulation of the 
critical values, while the power is comparable or even
slightly better than when using all intervals; this is presumably because the small loss
due to approximating the intervals is compensated by the smaller critical values due to
the fewer intervals that the statistic maximizes over.
In order to provide a fair comparison of the Bonferroni scan to the other four procedures
considered so far we simulated their critical values on the approximating set.
Figure~\ref{fig2} shows a comparison of
the critical values for sample size $n=10^6$. The plot shows that while the critical values
of {\sl BonferroniScan$_n$} are not as small as those of $SAC_n$ or the blocked scan, they
are still competitive, which is remarkable given that they were obtained with a Bonferroni
correction. We note that the critical values were simulated with $10^4$ simulation runs,
which is possible for the sample size $n=10^6$ because the approximating set has close to
$O(n)$ intervals rather than $O(n^2)$. The Bonferroni scan can be readily employed for
even much larger sample sizes since it does not require any simulations to find its critical values.

\section{Performance comparison} \label{sims}

This section compares the performance of the three calibrations with the traditional scan
in terms of the realized exponent $e_n(|I_n|)$. As explained in Section~\ref{optimality},
the realized exponent is a standardized measure of the smallest change in mean that the 
calibration can reliably detect in the model (\ref{model}).
Therefore, a calibration that performs well for this task should have a small realized
exponent, and ideally $e_n(|I_n|)$ should be close to the bound 1
given by the asymptotic detection boundary (\ref{consist}) for all signal lengths $|I_n|$.
We investigate this both theoretically and experimentally for the calibrations considered
so far.

\begin{Theorem} \label{asopt}
Each of the calibrations $DS_n$, $DS_n^{app}$, $SAC_n$, $SAC_n^{app}$, {\sl BlockedScan$_n$},
{\sl BlockedScan$_n^{app}$} and the Bonferroni scan has a realized exponent $e_n(|I_n|)$
that satisfies
$$
e_n(|I_n|) \ \leq \ 1+ \frac{b}{\sqrt{ \log \frac{n}{|I_n|}}}
$$
for all $|I_n| \in [1,n^p]$, $p \in (0,1)$, where $b=b(p)$ is a constant.
\end{Theorem}
 
We note that the traditional scan
$Scan_n$ does not satisfy such a bound as can be seen from the results in Chan and Walther~(2013).
Theorem~\ref{asopt} continues to hold if the elevated mean in model (\ref{model}) is not constant.
Then one simply needs to replace $\mu_n$ by the average of the means on $I$ in the definition of 
$e_n(|I_n|)$. This is readily seen from the proof of Theorem~\ref{asopt} as the test statistic is
a linear function of $\sum_{i \in I} Y_i$.
Furthermore, as explained in the proof, it is possible to sharpen the result of Theorem~\ref{asopt}
by deriving the appropriate value of $b$ for each calibration. Such a result would provide a theoretical
explanation of why $DS_n$ has an inferior performance for small $|I_n|$
and therefore complement existing optimality results that are not sensitive enough to discern
this effect, as was explained in Section~\ref{optimality}. We did not pursue this further since
our focus is the exact evaluation of $e_n(|I_n|)$ given below. 

In order to analyze the finite sample performance of these calibrations we evaluated
their realized exponents with simulations. We evaluated all calibrations on the approximating set,
for the reasons explained at the end of Section~\ref{methods}.
We first computed the realized exponent $e_n(|I_n|)$ for sample size $n=10^4$ using $10^4$ simulations.
Table~\ref{table1} shows $e_n(|I_n|)$ for a representative selection
of interval lengths $|I_n|$ for the various calibrations discussed in the previous section.
The column with signal length $I_n=1$ confirms the
poor performance of $DS_n$ discussed in the Introduction: the realized exponent is 1.80,
while that of the traditional scan is 1.41. In contrast, the performance of $SAC_n$ and of the 
blocked scan is only slightly inferior for the smallest scales while they increasingly outperform
the traditional scan from signal length 15 onwards. The Bonferroni scan  does not perform quite as
well as the $SAC_n$ and the blocked scan, but it is still competitive, providing a clear improvement
over the traditional scan for $|I_n| \geq 50$ without incurring the large penalty at small scales
that $DS_n$ does. We feel that this makes the Bonferroni scan an attractive choice since unlike
all the other calibrations, it does not require an ancillary approximation to obtain its critical 
values and it is therefore particularly simple to use.

We computed the realized exponents also for sample size $n=10^6$ and arrived at the same
conclusions, see Table~\ref{table2} and the discussion in Section~\ref{optimality}. 

Comparing the table for sample size $n=10^4$ with that for $n=10^6$ shows that
$e_n(|I_n|)$ converges very slowly to its asymptotic bound of 1, with many values
of $e_n(|I_n|)$ being closer to 2 rather than 1 even when $n=10^6$. This shows the merit
of evaluating the performance with the finite sample criterion $e_n(|I_n|)$ rather than 
by establishing an asymptotic
result as is typically done in the literature.
\smallskip

{\sl Details about the simulation study:} 
For each simulation of model (\ref{model}), the $Z_i$ were simulated and for a given length $|I_n|$
the start point of $I_n$ was chosen at random. Then a bisection search was used to find the smallest
$\mu_n$ for which the test rejects. The 80th percentile of these values is the Monte Carlo
estimate of $\mu_{min}(n,|I_n|)$, from which $e_n(|I_n|)$ is obtained from (\ref{realizedexp}).
We found that the conclusions of the simulation study are not sensitive to the chosen significance level
or the choice of 80\% power.

\begin{table}  
\centering
\begin{tabular}{l|cccccccc}
Signal length $|I_n|$ & 1 & 5 & 10 & 15 & 50 & 100 & 500 & 1000 \\ \hline
$Scan_n^{app}$ & 1.41 &1.58 &1.74 &1.85 &2.19 &2.47 &3.48 &4.34\\
$DS_n^{app}$   & 1.80 &1.79 &1.86 &1.90 &1.92 &1.94 &2.08 &2.25\\
$SAC_n^{app}$  & 1.41 &1.62 &1.76 &1.85 &2.04 &2.18 &2.73 &3.18\\
$BlockedScan_n^{app}$ & 1.49 &1.67 &1.83 &1.87 &1.91 &2.01 &2.43 &2.80\\
Bonferroni scan& 1.60 &1.81 &1.98 &2.03 &2.13 &2.25 &2.75 &3.17\\
\end{tabular}
\caption{The realized exponents $e_n(|I_n|)$ for various calibrations and sample size $n=10^4$.}
\label{table1}
\end{table}

\begin{table} 
\centering
\begin{tabular}{l|ccccccccccc}
Signal length $|I_n|$ & 1 & 5 & 10 & 20 & 50 & 100 & 500 & 1000 & 5000 & $10^4$ & $10^5$\\ \hline
$Scan_n^{app}$ & 1.32 &1.45 &1.55 &1.67 &1.79 &1.90 &2.27 &2.45 &3.20 &3.60 &6.26\\
$DS_n^{app}$   & 1.67 &1.69 &1.75 &1.79 &1.79 &1.78 &1.85 &1.87 &1.98 &2.04 &2.42\\
$SAC_n^{app}$  & 1.34 &1.50 &1.58 &1.67 &1.73 &1.77 &1.95 &2.03 &2.36 &2.54 &3.66\\
$BlockedScan_n^{app}$  & 1.40 &1.53 &1.63 &1.61 &1.71 &1.75 &1.90 &1.97 &2.22 &2.35 &3.11\\
Bonferroni scan& 1.46 &1.61 &1.72 &1.72 &1.83 &1.87 &2.05 &2.11 &2.45 &2.58 &3.59\\
\end{tabular}
\caption{The realized exponents $e_n(|I_n|)$ for various calibrations and sample size $n=10^6$.}
\label{table2}
\end{table}

\section{Tests against an epidemic alternative }
\label{epidemic}

As stated in the Introduction, the motivation for analyzing the abstract
Gaussian noise model (\ref{model}) is that key results in that abstract model typically carry over
to related and more realistic settings, such as when
the baseline is unknown, the noise variance is unknown,
or the error distribution is not normal. This section will discuss how the methodology introduced
in the previous sections applies in those settings. A key idea is that the statistic $T_I(\Yn)$
in model (\ref{model}) corresponds to the square root of the log likelihood ratio statistic, which
has a sub-Gaussian distribution under the null hypothesis, i.e. $f_n=0$ in (\ref{model}).

(\ref{model}) provides the conceptual framework for the problem of 
testing against a so-called epidemic (or square wave) 
alternative, see e.g. Levin and Kline~(1985), Siegmund~(1986) or Yao~(1993a,1993b):

One observes independent random variables $Y_i$ having mean $\mu_i$, $i=1\ldots,n$. Under the null
hypothesis all $\mu_i$ are equal:
\be \label{H0}
H_0:\ \mu_1=\ldots =\mu_n = \mu
\ee
while under the alternative there is an interval where the means are shifted:
$$
H_1:\ \mbox{for some interval $I \subset (0,n]$ and $\delta > 0$: } \mu_i = \begin{cases}
\mu +\delta &\mbox{if } i \in I\\
\mu      &\mbox{if } i \not\in I
\end{cases}
$$
where $\mu$ and $\delta$ are unknown nuisance parameters. For simplicity we will focus on the one-sided
alternative where $\delta >0$. This model is of interest because represents a template for an epidemic
state during the interval $I$ after which the normal state is restored, and it has also been employed in
other change-point problems such as the analysis of DNA copy number data, see Olshen et al.~(2004),
and in various problems in econometrics, see e.g. Broemeling and Tsurumi~(1987).

A number of different statistics have been proposed to analyze the above problem, see Yao~(1993a)
for a review and comparison. In particular, it was pointed out in Siegmund~(1986) and Yao~(1993a)
that there is a trade-off between these various statistics, with some statistics assigning more power
to smaller intervals $I$ at the cost of having less power at larger intervals, and it was proposed to
 impose a lower bound $|I| \geq \gamma n$ which gives the statistician the flexibility to trade power
between the scales of $|I|$. But this is precisely the problem that was addressed in the previous
sections, where it was shown how the use of scale-dependent critical values allows to achieve essentially
optimal inference for all scales $|I|$ simultaneously. 

Recall that the requirements for using that
methodology are as follows: The Sharpnack-Arias-Castro calibration assumes that the test statistic
$T_I(\Yn)$ is standard normal (or more generally sub-Gaussian with scale factor 1) under the null hypothesis.
A violation of this assumption will not invalidate the significance level or p-values, but it will lead
to a sub-optimal distribution of power across scales. As will be seen below, this assumption can be
readily met in many parametric and nonparametric settings, usually by employing the square-root of
the log likelihood ratio statistic as proposed by Rivera and Walther~(2013). The blocked scan does
not have such an assumption since it calibrates the significance level rather than the critical
value, but both the Sharpnack-Arias-Castro calibration and the blocked scan require that one can
simulate $\Yn$ under the null hypothesis in order to produce the critical values for each method.
The Bonferroni scan does not require this simulation, but it requires the knowledge of the distribution of
$T_I(\Yn)$ under the null hypthesis, or at least good bounds for tail probabilities, in order to
apply the Bonferroni correction.

\subsection{Normal observations with unknown baseline and known variance} 
\label{normal1}

We observe $Y_i \sim \mbox{N$(\mu_i,\sigma^2)$}$, $i=1,\ldots,n$, where $\sigma^2$ is known but the
baseline $\mu$ in (\ref{H0}) is unknown. This is the setting that is most commonly considered for the epidemic
model described above. The generalized likelihood ratio test rejects for large values of
\begin{align}
&\max_I \frac{1}{\sigma} \left( \sqrt{\frac{n -|I|}{n|I|}} \sum_{i \in I} Y_i -
\sqrt{\frac{|I|}{n(n-|I|)}} \sum_{i \not\in I} Y_i\right)  \nn \\
= &\max_I \ T_I(\Yn) \qquad \mbox{where }\ 
T_I(\Yn):=\ \frac{1}{\sigma} \left( \overline{Y}_I -\overline{Y} \right) \sqrt{\frac{n |I|}{n-|I|}}
\label{normalunk}
\end{align}
with $\overline{Y}_I:=\frac{1}{|I|}\sum_{i \in I}Y_i$ and $\overline{Y}:=\frac{1}{n}\sum_{i \in I}Y_i$;
see Hogan and Siegmund~(1986), Yao~(1993a), or Example~1 in Yu~(2020) for a step-by-step computation. It
is readily checked that $T_I(\Yn) \sim \mbox{N$(0,1)$}$ under $H_0$, so the generalized likelihood ratio test 
is given by the scan statistic calibration (\ref{scan}). Note that this setting is almost identical to the abstract
Gaussian model (\ref{model}) except that for disjoint $I_1,I_2$ the $T_{I_1}(\Yn)$ and $T_{I_2}(\Yn)$
are not quite independent because both contain the overall mean $\overline{Y}$. Therefore all three calibrations
described in Section~\ref{methods} are applicable. Indeed, simulations show that the realized exponents $e_n(|I_n|)$
of the various calibrations are very similar to those in Tables~\ref{table1} and \ref{table2} for the model
with known baseline, except for an increase (of roughly equal size for all calibrations) at the largest windows,
due to the aforementioned dependence.
Since the the null distribution of the statistic in (\ref{normalunk}) does not depend on the unknown
baseline $\mu$, the critical values for the Sharpnack-Arias-Castro calibration and for the blocked scan
can be computed by simulating the $Y_i$  from N$(0,\sigma^2)$.

\subsection{Normal observations with unknown baseline and unknown variance}

If $\sigma$ is unknown in the previous setting, then the generalized likelihood ratio test rejects for large values of
\be \label{normalunk2}
\max_I \ T_I(\Yn) \qquad \mbox{where }\
T_I(\Yn):=\ \frac{\overline{Y}_I -\overline{Y}}{\sqrt{\frac{1}{n-1} \sum_{i=1}^n(Y_i-\overline{Y})^2}}  
\sqrt{\frac{n |I|}{n-|I|}}
\ee
see Yao~(1993b), so $T_I(\Yn)$ is a studentized version of the statistic used in the previous setting.
At first sight, this would seem to be make it unadvisable to use the Sharpnack-Arias-Castro calibration
since the $t$-statistic has heavy algebraic tails rather than the expontial sub-Gaussian tails required
for that calibration. However, $\overline{Y}_I$ is not just studentized by also empirically centered
with the overall mean $\overline{Y}$ and it turns out that, perhaps surprisingly, this restores the
normal tail bound: It is shown in Walther~(2021) that
$$
\Pr (T_I(\Yn) >t)\ \leq \ \Pr (\mbox{N$(0,1) >t)$ \quad for } 
\begin{cases}
t\geq 2.5 \mbox{ and } n\geq 10, \ \mbox{ or} \\
t\geq 2.75 \mbox{ and } n\geq 6.
\end{cases}
$$
This makes it again possible to apply all three calibrations as in the abstract Gaussian model (\ref{model}).
The critical values for the Sharpnack-Arias-Castro calibration and for the blocked scan
can be computed by simulating the $Y_i$  from N$(0,1)$ since the null distribution of the test statistic
does not depend on the unknown baseline $\mu$ and $\sigma$.

\subsection{Observations from a natural exponential family}

Let $Y_i\sim f_{\theta_i}$, $i=1,\ldots,n$, be independent observations from a regular one-dimensional natural 
exponential family $\{f_{\theta}, \theta \in \Theta\}$, i.e. $f_{\theta}$ has a density with repect
to some $\sigma$-finite measure $\nu$ which is of the form $f_{\theta}(x)= \exp (\theta x
-A(\theta))\,h(x)$ and the natural parameter space $\Theta =\{\theta \in \R:\,
\int \exp ( \theta x) h(x) \nu (dx) < \infty \}$ is open.

The mean value parametrization $\mu_i=A'(\th_i)$ establishes the link to testing (\ref{H0}) against
an epidemic alternative. Note that the Gaussian detection problems
(\ref{model}) and \ref{normal1} are special cases of this model.
Exponential families are a natural and popular model for scanning problems, see e.g. Arias-Castro
et al.~(2011), Frick et al.~(2014) or K\"{o}nig et al.~(2020). One approach taken there is to
approximate the sum of observations by a normal distribution in order to use methodology applicable
for a normal model. However, a normal approximation does not imply a sub-Gaussian tail bound and therefore 
it is not clear how effective the Sharpnack-Arias-Castro calibration will be, especially if $|I|$ is small 
and so there are only few terms in the sum. For this reason we advocate instead for using 
the square root of twice the log likelihood ratio as test statistic $T_I(\Yn)$, as proposed by
Rivera and Walther~(2013). The advantage of this statistic is that it produces {\sl finite sample}
sub-Gaussian tail bounds, see Rivera and Walther~(2013) for the binomial case
and Walther~(2021) for the one-dimensional exponential family. The heuristic for this result is the  
well known fact that twice the log likelihood ratio is asymptotically chi-square distributed, hence the
square root will have approximately normal tails.

In more detail, the log likelihood ratio statistic for testing $H_0:\, \theta =\theta_0$ on $I$ is
\begin{align}
\logLRI (\theta_0) &:= \log \frac{\sup_{\theta \in \Theta} \prod_{i \in I} f_{\theta} (Y_i)}{
\prod_{i \in I} f_{\theta_0} (Y_i)} \nn \\
&= \sup_{\theta \in \Theta} \left( (\theta-\th_0) \sum_{i \in I} Y_i -|I| 
\Bl(A(\theta)-A(\th_0)\Br)\right) \label{logLRI}
\end{align}

See e.g. Barndorff-Nielsen~(1986) for inference based on $\sqrt{2 \logLRI (\theta_0)}$ or its signed
version for the case of one sided-alternatives. A normal approximation to this statistic is given in
Frick et al.~(2014) and K\"{o}nig et al.~(2020). Recent results show that a normal approximation
is not required for our purposes here since in fact a finite sample sub-Gaussian tail bound holds:
$\Pr_{\th_0} \left( \sqrt{2\, \logLRI (\th_0)} > x\right) \leq  2 \exp\left( -\frac{x^2}{2}\right)$
for $x>0$, see Walther~(2021).

In order to see the connection to the Gaussian model (\ref{model}), note that in the case where $f_{\th}$ 
is N$(\th,\sigma)$ with known $\th_0=0$ and $\sigma =1$, 
then (the signed) $\sqrt{2 \logLRI (\theta_0)}$ equals the statistic $T_I(\Yn)$ in (\ref{T}). 

Usually $\th_0$ is not known. Then the standard likelihood-based approach to test against a change of $\th$ on a 
given $I$ is the generalized likelihood ratio test, which divides the maximized likelihood under the alternative by
the maximized likelihood under the null hypothesis of a constant $\th$:

\begin{align}
{\rm logLR}&_{I,I^c}(\hat{\th}_{[n]}) \ :=\ \log \frac{\Bl(\sup_{\theta \in \Theta} \prod_{i \in I} f_{\theta} (X_i)\Br)
\left(\sup_{\theta \in \Theta} \prod_{i \not\in I} f_{\theta} (X_i)\right)}{
\sup_{\theta \in \Theta} \prod_{i=1}^n f_{\theta} (X_i)} \label{logLRIIc} \\
&= \sup_{\theta \in \Theta} \left( \theta \sum_{i \in I} X_i -|I| A(\theta)\right) +
\sup_{\theta \in \Theta} \left( \theta \sum_{i \not\in I} X_i -\Bl(n-|I|\Br) A(\theta)\right)-
\sup_{\theta \in \Theta} \left( \theta \sum_{i=1}^n X_i -n A(\theta)\right) \nn
\end{align}

We propose to use as test statistic $T_I(\Yn)=\sqrt{2 \logLRIIc}$. Importantly, this statistic
also satisfies a finite sample sub-Gaussian tail bound:
$\Pr_{\th_0} \left( \sqrt{2\, \logLRIIc } > x\right) \leq (4+2e) \exp\left( -\frac{x^2}{2}\right)$
for $x>0$, see Walther~(2021). (In the case of one-sided alternatives the corresponding
tail bounds for the signed square root are obtained by dividing the above bounds by 2.)
Therefore, the Sharpnack-Arias-Castro calibration will be effective for this statistic, as will be the
blocked scan (which does not require a sub-Gaussian tail bound). 
Critical values for these two calibrations can be computed by simulating from $f_{\hat{\th}_{[n]}}$,
where $\hat{\th}_{[n]}$ is the MLE under $H_0$, or by using the permutation distribution as proposed in Walther~(2010)
for the blocked scan, see also Section~\ref{permutation}. The Bonferroni scan requires the knowledge of
tail probabilities of $T_I(\Yn)$, which are available via the above tail bounds or by assuming that
a normal approximation holds (as is supported by simulations.)
 Simulations show that the former tail bounds are conservative by a factor of the order 10. Since the Bonferroni scan
will inherit this conservatism, one may prefer to use instead the standard normal approximation, provided one
is willing to work with an approximation instead of finite sample guarantees.

By way of illustration, if $f_{\th}$ is N$(\th,\sigma)$ with unknown baseline $\th_0$ and known $\sigma$,
then $T_I(\Yn)=\sqrt{2 \,\logLRIIc}$ equals the statistic $T_I(\Yn)$ given in (\ref{normalunk}).
As another example, the {\sl Bernoulli model} posits that the distribution of $Y_i$ is Bernoulli with unknown 
baseline probability $p$. The natural parameter for the exponential family is $\th=\log \frac{p}{1-p}$.
One computes
\be \label{logLR}
\logLRIIc\ = |I| \left( \overline{Y}_I \log \frac{\overline{Y}_I}{\overline{Y}} +(1-\overline{Y}_I)
\log \frac{1-\overline{Y}_I}{1-\overline{Y}}\right) +(n-|I|) \left( \overline{Y}_{I^c} \log
\frac{\overline{Y}_{I^c}}{\overline{Y}} +(1-\overline{Y}_{I^c}) \log
\frac{1-\overline{Y}_{I^c}}{1-\overline{Y}}\right)
\ee
where as before $\overline{Y}_I:=\frac{1}{|I|}\sum_{i \in I}Y_i$ and $\overline{Y}:=\frac{1}{n}\sum_{i \in I}Y_i$.
This statistic was proposed as a scan statistic by Kulldorff~(1997) and, despite its cumbersome form,
has been widely adopted for scanning problems in the computer science literature, see e.g. 
Neill and Moore~(2004a,2004b).

As a third example, we evaluated the realized exponents for the various calibrations 
in the case of a Poisson model where $Y_i \sim \mbox{Pois}(\mu_i)$ with an unkown baseline $\mu$.
The natural parameter is $\th=\log \mu$. One computes
$$
\logLRIIc\ = |I|\overline{Y}_I \bigl( \log \overline{Y}_I -1\bigr) +(n-|I|) \overline{Y}_{I^c} \bigl(
\log \overline{Y}_{I^c} -1\bigr) -n \overline{Y} \bigl( \log \overline{Y} -1 \bigr)
$$
Tables~\ref{table3} and \ref{table4} give the realized exponents for sample sizes $n=10^4$ and
$n=10^6$. The simulations were performed as described in Section~\ref{sims} using the baseline 
$\mu=1$, which was treated as unknown in the analysis. The realized exponents in these tables confirm
the conclusions about the various calibrations in the abstract Gaussian model in Section~\ref{sims}.

\begin{table}  
\centering
\begin{tabular}{l|cccccccc}
Signal length $|I_n|$ & 1 & 5 & 10 & 15 & 50 & 100 & 500 & 1000 \\ \hline
$Scan_n^{app}$ & 2.98 &2.19 &2.16 &2.21 &2.39 &2.64 &3.67 &4.80\\
$DS_n^{app}$   & 3.61 &2.63 &2.32 &2.16 &1.90 &1.90 &2.15 &2.44\\
$SAC_n^{app}$  & 2.94 &2.05 &1.96 &1.95 &1.98 &2.18 &2.68 &3.33\\
$BlockedScan_n^{app}$ & 2.99 &2.19 &2.13 &2.07 &1.89 &2.01 &2.53 &2.88\\
Bonferroni scan& 3.58 &2.71 &2.64 &2.56 &2.41 &2.49 &2.90 &3.63\\
\end{tabular}
\caption{The realized exponents $e_n(|I_n|)$ for the Poisson model and sample size $n=10^4$.}
\label{table3}
\end{table}

\begin{table} 
\centering
\begin{tabular}{l|ccccccccccc}
Signal length $|I_n|$ & 1 & 5 & 10 & 20 & 50 & 100 & 500 & 1000 & 5000 & $10^4$ & $10^5$\\ \hline
$Scan_n^{app}$ & 3.26 &2.24 &2.25 &2.21 &2.19 &2.22 &2.44 &2.59 &3.26 &3.67 &7.20\\
$DS_n^{app}$   & 4.61 &2.89 &2.53 &2.21 &2.02 &1.95 &1.88 &1.88 &1.89 &2.00 &2.68\\
$SAC_n^{app}$  & 3.30 &2.20 &2.09 &2.04 &1.96 &1.92 &1.95 &1.96 &2.37 &2.52 &4.05\\
$BlockedScan_n^{app}$  & 3.29 &2.45 &2.20 &2.00 &1.93 &1.94 &1.94 &1.93 &2.24 &2.41 &3.59\\
Bonferroni scan& 3.77 &2.65 &2.56 &2.06 &2.03 &2.03 &2.07 &2.12 &2.55 &2.64 &4.10\\
\end{tabular}
\caption{The realized exponents $e_n(|I_n|)$ for the Poisson model and sample size $n=10^6$.}
\label{table4}
\end{table}

We will not formally check that the theoretical optimality results about the detection boundary
(\ref{consist}) and Theorem~\ref{asopt} hold for the exponential family setting. We refer instead to the heuristic
that the results for the abstract
Gaussian noise model (\ref{model}) can be expected to carry over to related situations where inference
is based on the likelihood ratio statistc, see e.g. Brown and Low~(1996) for general results.

\subsection{Heteroscedastic observations from a normal or symmetric distribution}

In this section we consider observations that may not be identically distributed under the null
hypothesis. In particular, we will address the heteroscedastic normal model where the 
$Y_i \sim \mbox{N}(\mu_i,\sigma_i^2)$
are independent with unknown $\mu_i,\sigma_i$ and the $\sigma_i$ are not necessarily identical, 
but the methodology is more generally  
applicable to independent $Y_i$ that are symmetric around some unknown $\mu_i$ (and thus possibly heteroscedastic). 
To this end, we propose
to use the technique of {\sl self-normalization}, see de la Pe\~{n}a et al.~(2009).
The key idea is that if the $Y_i$ are independent (but not necessarily identically distributed)
and symmetrically distributed around zero, then the self-normalized sum
$\sum_{i=1}^n Y_i / \sqrt{\sum_{i=1}^n Y_i^2 }$ will satisfy a sub-Gaussian tail bound, due
to Hoeffding's inequality applied conditional on the $Y_i$. Moreover, recent results about Rademacher
sequences by Pinelis~(2012) and Bentkus and Dzindzalieta~(2015) imply that a tail bound holds that is in fact
close to that of a standard normal. In order to appreciate this result, let e.g. $Y_1$ be Cauchy,
let $Y_2$ have a $t$-distribution with 5 degrees of freedom, and let $Y_3$ have a double exponential distribution.
Then after self-normalization the sum $Y_1+Y_2+Y_3$ will obey a guaranteed tail bound that is very close
to the standard normal bound. 
The intuition behind this perhaps surprising result is that a large contribution by an individual $Y_i$
to the sum will be dampened by the same large contribution in the studentizing term.

In order to make this result useful for our setting where the $Y_i$ are symmetric around an unknown $\mu_i$,
we need to extend the tail bound to the case where the center of symmetry is not zero.
This would seem to be a hopeless problem since the symmetry around zero is key for establishing
the tail bound. However, the situation here is different from the usual setting for self-normalized
processes in that there are additional observations outside the scanning window. The idea is to subtract off
these observations in order to  eliminate the unknown center of symmetry: For $i \in I$ we compute
$$
\wt{Y}_i\ :=\ Y_i -\frac{1}{\# J_i} \sum_{j \in J_i} Y_j
$$
where $\{J_i,\,i \in I\}$ is a partition of $\{1,\ldots,n\} \setminus I$. For simplicity we assume
$n=|I| p$ for some integer $p \geq 2$, which can always be arranged by discarding some data if
necessary. Then we can arrange each index set $J_i$ to contain $p-1$ indices. Writing $Y=(Y_1,\ldots,Y_n)^T$
and $\wt{Y}=(\wt{Y}_i,i \in I)^T$, one sees from the definition of
$\wt{Y}_i$ that $\wt{Y}={\bs A}_I Y$ for a certain matrix ${\bs A}_I$, hence the
self-normalized sum $\sum_{i\in I} \wt{Y}_i / \sqrt{\sum_{i\in I} \wt{Y}_i^2 }$ is readily computed as
$$
\frac{n}{n-|I|} \frac{\sum_{i \in I} (Y_i-\overline{Y})}{\sqrt{Y^T A_I^T A_I Y}}\ =:\ T_I(\Yn)
$$
Note that this self-normalized sum has the same form as the likelihood ratio statistic in the
normal model (\ref{normalunk2}), the difference being the studentizing term and the scaling.
Studentizing by $\sqrt{Y^T A_I^T A_I Y}$ turns out to be effective for dealing with arbitrary heteroscedasticity:
Based on the above results about Rademacher sequences
it is shown in Walther~(2021) that if for some real numbers $\mu_I \leq \mu_{I^c}$ we have $\mu_i=\mu_I$ for all
$i \in I$ and $\mu_i =\mu_{I^c}$ for all $i \not\in I$, then
\be  \label{bound318}
\Pr \left( T_I(\Yn) \geq t \right) \ \leq \ \min \Bl( 3.18, g(t)\Br) \,\Pr \Bl({\rm N(0,1)} > t\Br)
\ee
for all $t>0$, where $g(t):=1+\frac{14.11 \phi(t)}{(9+t^2) (1-\Phi(t))} \ra 1$ as $t \ra \infty$.
Hence the tail of $T_I(\Yn)$ is at most 3.18 times that of a standard normal (and even less if $t$ is large
enough to so that $g(t)<3.18$) in the null case when the means are constant on $I$ and on $I^c$ and the means
are not elevated on $I$ (i.e. smaller or equal
than on $I^c$). Importantly, this result continues to hold when the
means $\mu_i$ vary on $I$ and on $I^c$, provided that either the means $\Ex \wt{Y}_i$ don't vary much,
or the $\mu_i$ don't vary much for $i \in I$ and for $i \not\in I$, see Walther~(2021) for details.

While the existence of this normal tail bound implies that the Sharpnack-Arias-Castro calibration
will be effective, it is not clear how one would go about simulating its critical values since
the $\sigma_i$ are unknown. This also applies for the blocked scan. Therefore only the Bonferroni scan
will be readily applicable in this situation, with tail probabilities for the $T_I(\Yn)$ obtained via
the bound (\ref{bound318}). Simulations show that this tail bound is conservative for self-normalized
sums, but some conservatism is arguably a small price to pay in turn for being able to address
the difficult heteroscedastic case at all. For example, the recent work of Enikeeva~(2018) addresses the
Gaussian detection problem (\ref{model}) while allowing $\sigma$ to change together with $\mu$ on $I$,
but it is assumed that $\sigma$ is constant and known on both $I$ and on $I^c$.

Since the above approach appears to be new, it is of interest to ascertain its performance.
Set $\sigma_I^2:=\frac{1}{|I|}\sum_{i \in I} \sigma_i^2$.
It is shown in Walther~(2021) that if the $\sigma_i$ do not change quickly, namely $\sigma_j^2/\sigma_I^2 \leq
S \sqrt{\max_{i \in I} (j-i)}$ for all $j \in
\{1,\ldots,n\}$ and some $S>0$, then the
Bonferroni scan has asymptotic power 1 provided
$$
\frac{1}{|I|} \sum_{i \in I} \mu_i \,-\, \frac{1}{n-|I|} \sum_{i \not\in I} \mu_i\ \geq\ 
\sqrt{\frac{(2+\eps_n)\, \sigma_I^2\, R_I\, \log \frac{n}{|I|}}{|I|}}
$$ 
with $\eps_n \sqrt{\log \frac{n}{|I|}} \ra \infty$, where
$$
R_I\ :=\ \frac{\sum_{i \in I} {\rm Var} \wt{Y}_i}{\sum_{i \in I} \sigma_i^2}\ \leq\ 
1+2S \sqrt{\frac{|I|^2}{n}}
$$

So in the special case of constant $\sigma_i =:\sigma$ for all $i\in I$ and $|I| \leq \frac{n^{1/2}}{\log n}$, 
we get $\sigma_I=\sigma$ and $R_I=1+o(\eps_n)$ and hence, remarkably, the detection threshold of this more general
method is the same as the optimal detection boundary (\ref{consist}) in the homoscedastic case, even though the 
$\sigma_i$ may vary outside $I$ (subject to the above growth condition).

\subsection{Observations from an arbitrary distribution that are exchangeable under the null hypothesis:
Tests based on permutations, ranks or signs}
\label{permutation}

If one wishes to relax the distributional assumption on the $Y_i$ further, then a pertinent approach
is the simple and popular permutation test. This test produces exact finite sample significance statements due 
to the underlying invariance argument, see Hoeffding~(1952) or Romano~(1989). This invariance requires
that the joint distribution of the $Y_i$ is exchangeable (for example i.i.d.) under the null hypothesis, 
but it does not require an assumption on the distribution of the $Y_i$. Thus this setting is more general
than the first three settings considered above, but it does not include the previous one that allows heteroscedasticity.
It was shown in Arias-Castro et al.~(2018) that for observations from an exponential family the permutation
test will asymptotically detect epidemic alternatives if the parameter differs from $\th_0$ by at least
$\tau \sqrt{2 \log n}/\sqrt{|I|}$ with $\tau >1$. 
Thus for alternatives on very short intervals $|I|$ the permutation test will obtain the optimal
detection boundary (\ref{consist}), even without using the knowledge that the observations are from
an exponential family. However, if $|I|$ is not very short, then the optimal boundary (\ref{consist})
is smaller. As explained in Section~\ref{optimality}, if for example $|I|/n=n^{-1/2}$, then the exponent of effective
dimension can be reduced from 2 to 1, provided one employs scale-dependent critical values. It is therefore
of interest to see whether any of the scale-dependent calibrations is applicable to the permutation scan.

If one uses as test statistic the sum of the observations, then it turns out that there is no guarantee
that the Sharpnack-Arias-Castro calibration will provide an optimal trade-off between the scales,
and the Bonferroni scan cannot be applied since it is not clear how to obtain the required tail probabilities (in
a computationally efficient manner). Therefore the blocked scan will be the go-to method in this case, as
proposed by Walther~(2010). The reason why the Sharpnack-Arias-Castro calibration may be suboptimal is that
the permutation distribution of the sum satisfies only a Bernstein-type
tail bound rather than the tighter sub-Gaussian bound, see e.g. Bardenet and Maillard~(2015) or Lemma~2
in Arias-Castro et al.~(2018). Indeed, it is known that the permutation distribution of the sum is close
to its original distribution, see e.g. Romano~(1989). So if the $Y_i$ have a Cauchy distribution, then the
permutation distribution of their sum will be close to Cauchy, and the Cauchy tail is quite different
from  the normal tail for which the Sharpnack-Arias-Castro calibration is designed.

Alternatively, if one wishes to work with statistics based on ranks or signs, then the other two calibrations
also become readily applicable. Jung and Cho~(2015) and Arias-Castro et al.~(2018) propose to use
as test statistic the sum of the ranks in $I$ rather than the sum of the observations. Let
$\mathbf{R}_n=(R_1,\ldots,R_n)$ denote the ranks of $(Y_1,\ldots,Y_n)$. Jung and Cho~(2015) propose to 
scan with the Wilcoxon rank-sum statistic 
$T_I(\mathbf{R}_n)=\sqrt{\frac{12 |I|}{(n+1)(n-|I|)}} \left( \frac{1}{|I|} \sum_{i \in I} R_i -\frac{n+1}{2}\right)$.
The Wilcoxon rank-sum statistic is known to be close to normal even for relatively small sample sizes.
Indeed, it follows from Corollary~6.6 in D\"{u}mbgen~(2002) that $\Pr \left( \sqrt{\frac{n-|I|}{n+1}}\, T_I(\mathbf{R}_n)
>t\right) \leq \exp \left(-t^2/2\right)$, so the sub-Gaussian tail bound holds after taking out the correction factor
for sampling without replacement (which becomes negligible for smaller  $|I|$). Therefore the 
Sharpnack-Arias-Castro calibration will be effective for these rank statistics.
Likewise, the Bonferroni scan is applicable since the tail probabilities of $T_I(\mathbf{R}_n)$ can be
evaluated with a standard normal distribution or with the above sub-Gaussian tail bound.

Classifying the observations in $I$ according to whether they are above or below the overall median
results in a local sign test statistic $T_I(\Yn)=\sum_{i \in I} {\bs 1}(Y_i \geq$ median$(\Yn))$.
The sign test has the advantage that it is popular and easily interpretable by non-experts,
but comes at the cost of sacrificing some power. The permutation distribution of $T_I(\Yn)$
is hypergeometric, corresponding to drawing $|I|$ labels without replacement from $n$ of which half are 1 and half
are 0 if $n$ is even. Therefore the Bonferroni scan is readily applicable, as are the other two
calibrations. For the Sharpnack-Arias-Castro calibration it is helpful to use the following transformation
of the hypergeometric distribution in order to produce cleaner sub-Gaussian tails: Let $T_I(\Yn)=\sqrt{2 \,\logLRIIc}$
where $\logLRIIc$ is given in (\ref{logLR}) with $Y_i$ replaced by ${\bs 1}(Y_i \geq$ median$(\Yn))$.
Then it follows from Theorem~4 in Walther~(2010) that $T_I(\Yn)$ has the desired sub-Gaussian tail.

\section{Other settings}
\label{other}

The results of this paper were derived for the Gaussian sequence model (\ref{model})
because this allows to focus on the main ideas. However, the conclusions and methodology
can be carried over to other settings where scan statistics are used, beyond the various
regression settings discussed in Section~\ref{epidemic}.
For example, one case of particular interest in the
literature is the setting where one observes an (in)homogeneous Poisson Process and the problem
is to detect an interval where the intensity is elevated compared to a known baseline, i.e. one
is looking for an interval with an unusually large number of events, see Glaz et al. (2001).
Conditioning on the total number of observed events allows to eliminate certain nuisance parameters
and shows that the problem is equivalent to testing whether $n$ i.i.d. observations arise from
a known density $f_0$ (which w.l.o.g. can be taken to be the uniform density on $[0,1]$)
 versus the alternative that $f_0$ is elevated on an interval $I$:
\be  \label{density}
f_{r,I}(x)\ =\ \frac{r 1(x \in I) + 1(x \in I^c)}{rF_0(I) +F_0(I^c)}f_0(x),
\ee
so the problem becomes testing $r=1$ vs. $r>1$, see Loader (1991) and Rivera and Walther (2013).
The results in the latter paper suggest that the heuristics, methodology and optimality results
in the density/intensity model (\ref{density}) are quite analogous to the Gaussian sequence
model (\ref{model}). In particular,
it was pointed out in Rivera and Walther~(2013) that the {\sl square root} of the log likelihood
ratio statistic has a subgaussian tail, which allows to transfer the methodology from the Gaussian
sequence model just as in the exponential family case above. In the density setting the empirical measure 
plays the role that the interval length $|I_n|/n$ plays in the regression setting, and
the approximating set in Section~\ref{methods} will omit the first block as it is well known that at 
least $2 \log n$ observations are necessary for sensible inference in the density setting.

Some of the most important  applications of scan statistics concern the multivariate setting.
In that situation it is particularly important to evaluate the performance with a finite
sample criterion such as $e_n(|I_n|)$ because of the large-scale multiple testing problem
involved: While there are of the order $\sim n^2$
intervals that contain distinct subsets of $n$ points sampled from $U[0,1]$,
there are $\sim n^{2d}$ distinct axis-parallel rectangles
that contain distinct subsets of $n$ points sampled from $U[0,1]^d$. Nevertheless,
Arias-Castro et al.~(2005) show (for observations on a regular d-dimensional grid)
that for many relevant classes of scanning windows, such as axis-parallel rectangles
and balls, the effective dimension of the multiple testing problem is essentially  linear
in the sample size, i.e. the problem is not fundamentally more difficult than in the
univariate model (\ref{model}).
Moreover, it is shown in Walther~(2010) that 
if one employs scale-dependent
critical values (the blocked scan) rather than the traditional scan as in Arias-Castro et al.~(2005),
then it is possible to overcome the `curse of dimensionality':
if the signal is supported on a lower-dimensional marginal, then the d-dimensional blocked scan
has essentially the same asymptotic detection power as an optimal lower-dimensional  test
would have, so there is no fundamental penalty for scanning in the higher dimensional space.
The optimality results of Arias-Castro et al.~(2005) and Walther~(2010) are asymptotic.
Since the questionable utility of of asymptotic results as described in Section~\ref{optimality}
will be an even bigger concern in the multivariate setting,
it is of interest to reexamine the methodology with a finite
sample criterion such as $e_n(|I_n|)$. We note that all of the methodology developed in this
paper can be carried over to the multivariate setting. 
In fact, the approximating set in Section~\ref{methods}
is essentially the univariate version of the $d$-dimensional approximating set introduced in 
Walther~(2010).

\newpage

\appendix
{\LARGE\bf Appendix}

\section{Proof of Theorem~\ref{thm1}}  

For simplicity we write $\Yn(I):=T_I(\Yn)= 
\frac{\sum_{i \in I} Y_i}{\sqrt{|I|}}$ and likewise for $\Zn(I)$. Model (\ref{model}) gives
\be  \label{a}
\Yn(I)\ =\ \Zn(I) + \mu_n \frac{|I \cap I_n|}{\sqrt{|I|}}
\ee
Hence 
\begin{align}
\max_I \Yn(I) & = \max \Bl(\max_{I:I \cap I_n = \emptyset} \Yn(I), \max_{I:\frac{|I \cap I_n|}{
\sqrt{|I||I_n|}} \geq \frac{3}{4}} \Yn(I), \max_{I:0<\frac{|I \cap I_n|}{\sqrt{|I||I_n|}}
   < \frac{3}{4}} \Yn(I) \Br) \nn \\
& \leq \max \Bl( \max_I \Zn(I), \max_{I:\frac{|I \cap I_n|}{\sqrt{|I||I_n|}} \geq \frac{3}{4}} \Zn(I)
  +\mu_n \sqrt{|I_n|}, \max_{I:0<\frac{|I \cap I_n|}{\sqrt{|I||I_n|}} < \frac{3}{4}}\Zn(I)
   +\frac{3\mu_n \sqrt{|I_n|}}{4} \Br). \label{maxineq}
\end{align}

\begin{Lemma}  \label{LemmaA}
\be \label{b}
\max_{I:\frac{|I \cap I_n|}{\sqrt{|I||I_n|}} \geq \frac{3}{4}} \Zn(I)\ \stackrel{d}{\leq}\ R
\ee
where $R$ is a universal random variable whose support is the real line. Further
\be \label{c}
A_n\ :=\ \max_{I:I \cap I_n \neq \emptyset} \Zn(I) \ \mbox{ satisfies } A_n -\frac{1}{4} \sqrt{2 \log n}
\stackrel{p}{\ra} - \infty.
\ee
\end{Lemma}

\nin The lemma and (\ref{maxineq}) yield
$$
\Pr_{f_n} ( Scan_n(\Yn) > \kappa_n )  \leq \Pr (\max_I \Zn(I) > \kappa_n) + \Pr\Bl(R > \kappa_n
  - \mu_n \sqrt{|I_n|}\Br) + \Pr\Bl(A_n > \kappa_n - \frac{3 \mu_n \sqrt{|I_n|}}{4}\Br)
$$
Since $Scan_n(\Yn)$ is asymptotically powerful against $\{f_n\}$ the LHS 
converges to $1$ and the first term on the RHS converges to $0$. Hence, writing 
$b_n := \mu_n \sqrt{|I_n|} - \kappa_n$ and observing $\kappa_n \geq \sqrt{2 \log n}$
eventually by (2) in Kabluchko~(2011):
$$
1 \ \leq \ \liminf_n \Bl( \Pr(R > -b_n) + \Pr\Bl(A_n - \frac{1}{4} \sqrt{2 \log n} > -\frac{3}{4} 
 b_n\Br)\Br)
$$
This implies $b_n \ra \infty$ since the support of $R$ is the real line and 
$A_n -\frac{1}{4} \sqrt{2 \log n} \stackrel{p}{\ra} - \infty$.
But then we have for any sequence $\tilde{\kappa}_n = \kappa_n +O(1)$:
\begin{align}
\Pr_{f_n} (Scan_n(\Yn) > \tilde{\kappa}_n) & \geq  \Pr_{f_n}(\Yn(I_n) > \tilde{\kappa}_n) \\
& = \Pr\Bl( N(0,1) +\mu_n \sqrt{|I_n|} > \tilde{\kappa}_n \Br)\ \ \mbox{by (\ref{a})}\\
& \ra 1 \ \ \mbox{ since $\mu_n \sqrt{|I_n|} -\kappa_n \ra \infty$},
\end{align}
and $\Pr_0(Scan_n(\Yn) > \tilde{\kappa}_n) \leq \Pr(Scan_n(\Zn) > \kappa_n ) \ra 0$
if $\tilde{\kappa}_n \geq \kappa_n$. Hence $Scan_n$ is also asymptotically powerful against $\{f_n\}$
when employing the critical values $\tilde{\kappa}_n$. In fact, this is also true for any
sequence of critical values $\tilde{\kappa}_n =\sqrt{2\log n} +O(1)$ with $O(1)\geq 1 $ (say),
as can be seen by instead taking $b_n := \mu_n \sqrt{|I_n|} - \sqrt{2 \log n}$ above
and using (\ref{SV}).

We note that the above proof uses the fact that the type I error probability goes to 0. Some
related optimality results in the literature follow the approach of D\"{u}mbgen and Spokoiny~(2001)
and establish asymptotic power 1 at a fixed significance level. The main conclusion of this
theorem, namely that asymptotic optimality leaves a leeway of size $O(1)$ for the critical
value, will typically hold also in that setting and for related statistics, as can be seen
by the inspecting the proofs.

It remains to prove Lemma~\ref{LemmaA}. In order to prove (\ref{b}) we will show that
\be \label{d}
\max_{I:\frac{|I \cap I_n|}{\sqrt{|I||I_n|}} \geq \frac{3}{4}} \Zn(I) \ \stackrel{d}{\leq}\ 
R_1+R_2+R_3
\ee
where the $R_i$ are independent universal random variables and the support of $R_2$ is the
real line, hence the support of $R_1+R_2+R_3$ also equals the real line.

It is straightforward to check that the condition on $I$ implies $|I \cap I_n| \geq 
\frac{9}{16}|I_n|$ and $|I| \leq \frac{16}{9} |I_n|$. Therefore there exist intervals
$S_n$ and $L_n$, each having integers as endpoints and depending only on $I_n$, such that
$|S_n| \geq |I_n|/8$, $|L_n| \leq 5|I_n|$ and $S_n \subset I \subset L_n$ for every $I$.
(Let $S_n$ be the smallest such interval whose midpoint equals that of $I_n$, and construct
$L_n$ by moving each endpoint of $I_n$ outward by $2|I_n|$, then intersect the resulting
interval with $(0,n]$).
Hence we can write $I$ as the union of three disjoint intervals $I=I_{left}\cup S_n \cup I_{right}$,
where $I_{left}, I_{right}$ might be empty and $|I_{left}|,|I_{right}| \leq 2|I_n|$. So
$$
\Zn(I)\ =\ \frac{\sum_{i \in I_{left}} Z_i}{\sqrt{|I|}} +\frac{\sum_{i \in S_n} Z_i}{\sqrt{|I|}}
+\frac{\sum_{i \in I_{right}} Z_i}{\sqrt{|I|}}
$$
The three terms are independent. As for the middle term, $\frac{|S_n|}{|I|} \in
[\frac{1}{5\cdot 8},1]$ implies
\begin{align*}
\max_I \frac{\sum_{i \in S_n} Z_i}{\sqrt{|I|}} &= \max_I \sqrt{\frac{|S_n|}{|I|}} \Zn(S_n)\\
& \leq \sqrt{\frac{1}{40}} \Zn(S_n) \,{\bs 1} (\Zn(S_n) <0) + \Zn(S_n) \,{\bs 1}(\Zn(S_n) \geq 0) \\
& \stackrel{d}{=} \sqrt{\frac{1}{40}} Z \,{\bs 1} (Z <0) + Z \,{\bs 1}(Z \geq 0)\ =:\ R_2 
\end{align*}
where $Z \sim N(0,1)$. Clearly the support of $R_2$ is real line. As for the third term,
$\sum_{i \in I_{right}} Z_i$ is the sum of the first $|I_{right}|$ $Z_i$s following the 
right endpoint of $S_n$. Using $|I_{right}| \leq 2 |I_n|$ and $|I| \in [|I_n|/8, 5|I_n|]$
we get
\begin{align*}
\max_I \frac{\sum_{i \in I_{right}} Z_i}{\sqrt{|I|}} & \stackrel{d}{\leq} \max_{k=0,1,\ldots,2|I_n|}
  \max_{f_k \in [\frac{1}{8} |I_n|,5|I_n|]} \frac{\sum_{i=1}^k Z_i}{\sqrt{f_k}} \\
& \stackrel{d}{\leq} \sup_{0 \leq t \leq 2|I_n|} \sup_{f_t \in [\frac{1}{8} |I_n|,5|I_n|]}
  \frac{W(t)}{\sqrt{f_t}}\ \ \ \mbox{ where $W(t)$ is Brownian motion}\\
& = \sqrt{16} \sup_{0 \leq t \leq 2|I_n|} \frac{W(t)}{\sqrt{2|I_n|}} \ \ \ 
  \mbox{ since the sup is $\geq 0$ a.s.}\\
& =: \sqrt{16} R_3
\end{align*}
The distribution of $R_3$ is known to be that of a standard normal conditional on being positive.
The first term is analogous, proving (\ref{d}) since stochastic ordering is preserved when taking
sums of independent random variables.

As for (\ref{c}), write $I_n=(j_n,k_n]$ and note that $I \cap I_n \neq \emptyset$ implies that
$I \subset I_n$ or the left endpoint $j_n \in I$ or the right endpoint $k_n \in I$. Hence
\be  \label{C*}
\max_{I:I \cap I_n \neq \emptyset} \Zn(I) \leq \max \Bl(\max_{I \subset I_n} \Zn(I), 
  \max_{I:j_n \in I} \Zn(I), \max_{I: k_n \in I} \Zn(I)\Br)
\ee
By Lemma~1 of Chan and Walther~(2013), $\max_{I \subset I_n} \Zn(I) \stackrel{d}{\leq} L +
\sqrt{2 \log (e|I_n|)}$ for some universal random variable $L$. So if $|I_n| \leq n^p$ with
$p<\frac{1}{16}$, then $\max_{I \subset I_n} \Zn(I) - \frac{1}{4} \sqrt{2 \log n}
\stackrel{p}{\ra} - \infty$. Next,
\begin{align*}
\max_{I:j_n \in I} \Zn(I) & = \max_{a \in \{0,\ldots,j_n-1\}} \max_{b \in \{j_n+1,\ldots,n\}}
\frac{\sum_{i=a+1}^{j_n} Z_i + \sum_{i=j_n+1}^b Z_i}{\sqrt{b-a}} \\
& \leq \max_{a \in \{0,\ldots,j_n-1\}} \frac{\Bl|\sum_{i=a+1}^{j_n} Z_i|}{\sqrt{j_n-a}}  + 
  \max_{b \in \{j_n+1,\ldots,n\}} \frac{\Bl|\sum_{i=j_n+1}^b Z_i\Br|}{\sqrt{b-j_n}}
\end{align*}
These two terms are independent since they involve disjoint sets of $Z_is$. Each term is
$\stackrel{d}{\leq} \max_{1 \leq k \leq n} \frac{\sum_{i=1}^k Z_i}{\sqrt{k}}=:L_n$, hence
$\max_{I:j_n \in I} \Zn(I) \stackrel{d}{\leq} L_n +L_n'$, where $L_n$ and $L_n'$ are i.i.d.
copies, since stochastic ordering is preserved when taking sums of independent random variables.
By a theorem of Darling and Erd\"{o}s~(1956), $L_n =O_p(\sqrt{\log \log n})$. The third term
in (\ref{C*}) is bounded analogously. (\ref{c}) follows from (\ref{C*}) and the above bounds. $ \Box$

\section{Proof of Theorem~\ref{asopt}}

We will first establish the claim for
$DS_n$, $SAC_n$ and {\sl BlockedScan$_n$}.
For any signal $f_n$ from model (\ref{model}) with support $I_n=(j_n,k_n]$ we have
$T_{j_nk_n}(\Yn)  = \sqrt{|I_n|} \,\mu_n +T_{j_nk_n}(\Zn)$
by (\ref{a}). So if a generic test assigns critical values $c_{jkn}(\alp)$ to the $T_{jk}(\Yn)$,
then 
\begin{align}
\Pr_{f_n} (\mbox{test rejects}) & \geq \Pr_{f_n}\Bl( T_{j_nk_n}(\Yn) >c_{j_nk_nn}(\alp)\Br) \nn \\
& = 1- \Phi \Bl(c_{j_nk_nn}(\alp) - \sqrt{|I_n|} \,\mu_n \Br)\ \ 
  \mbox{ since $T_{j_nk_n}(\Zn) \sim N(0,1)$}  \label{*A}\\
& \geq 80\% \ \ \ \mbox{ provided $\sqrt{|I_n|} \,\mu_n \geq c_{j_nk_nn}(\alp)+z_{0.2}$},\nn
\end{align}
where $z_{0.2}$ denotes the 80th percentile of $N(0,1)$.
Since this
inequality holds uniformly for such $f_n$ it follows that the smallest detectable value 
$\sqrt{|I_n|} \,\mu_{min}(n,|I_n|)$ for this generic test is not larger than $c_{j_nk_nn}(\alp)+z_{0.2}$.
So if the $c_{jkn}(\alp)$ satisfy
\be  \label{*R}
c_{jkn}(\alp) - \sqrt{2 \log \frac{en}{k-j}} \ \leq\ b \ \ \mbox{ for all $0 \leq j<k \leq n$
with $k-j \leq n^p$}
\ee 
for some number $b$, then
$$
e_n(|I_n|)\ =\ \frac{\left(\sqrt{|I_n|} \,\mu_{min}(n,|I_n|)\right)^2 }{2 
\log \frac{en}{|I_n|}} \ \leq \ 1+\frac{2(b+z_{0.2})
}{\sqrt{2 \log \frac{en}{|I_n|}}} +\frac{(b+z_{0.2})^2}{2 \log \frac{en}{|I_n|}} \
 \leq \ 1+\frac{\sqrt{2}(b+z_{0.2}) +\frac{1}{2}(b+z_{0.2})^2}{
\sqrt{ \log\frac{en}{|I_n|}}}.
$$

Hence the claim of the theorem follows for a particular calibration $c_{jkn}(\alp)$
once (\ref{*R}) is established for that calibration. We will now check this condition for
the calibrations used by $DS_n$, $SAC_n$ and {\sl BlockedScan$_n$}. It follows from
Thm.~2.1 in D\"{u}mbgen and Spokoiny~(2001) that $DS_n(\Zn) =O_p(1)$. Hence its critical
value for $T_{jk}$ is given by $c_{jkn}(\alp)=\sqrt{2 \log \frac{en}{k-j}}+\kappa_n(\alp)$
with $\kappa_n(\alp)=O(1)$, so (\ref{*R}) holds. Next we check this condition for $SAC_n$.
Since the penalty term in $SAC_n$ is larger than
that in $DS_n$ we obtain $SAC_n(\Zn)=O_p(1)$ and hence the $(1-\alp)$-quantile of 
$SAC_n(\Zn)$ also stays bounded in $n$. So (\ref{*R}) holds for $SAC_n$ provided that
$$
\sqrt{2\log \Bl[\frac{en}{k-j}(1+\log (k-j))^2\Br]} \ \leq \ \sqrt{2 \log \frac{en}{k-j}}
+O(1) \ \ \mbox{ for all $0 \leq j<k \leq n$ with $k-j \leq n^p$}.
$$
But this is a consequence of $k-j \leq n^p$ since using $\sqrt{1+x}\leq 1+\frac{1}{2}x$ one
sees that the first square root is not larger than
$$
\sqrt{2 \log \frac{en}{k-j}} +\sqrt{2}\frac{\log \log (e(k-j))}{\sqrt{\log \frac{en}{k-j}}}
\ \leq \ \sqrt{2 \log \frac{en}{k-j}}  +\sqrt{\frac{2}{1-p}}.
$$
As an aside, (\ref{*A}) also shows how much $SAC_n$ loses on the largest scales
$k-j \approx n$ when compared to $DS_n$: Then the DS penalty is $\sqrt{2 \log \frac{en}{k_n-j_n}}
\approx \sqrt{2}$, while the SAC penalty is approximately $\sqrt{4 \log \log n}$.
So in order to achieve 80\% power via the condition $\sqrt{|I_n|} \,\mu_n \geq c_{j_nk_nn}(\alp)
+z_{0.2}$ in (\ref{*A}), it follows that some constant $\sqrt{|I_n|} \,\mu_n$ is sufficient for
$DS_n$ while $SAC_n$ needs $\sqrt{|I_n|} \,\mu_n$ to grow with $n$ but only
at the very slow rate $\sqrt{4 \log \log n}$.

Furthermore, proceeding similarly as in Chan and Lai~(2006) and Kabluchko~(2011)
 it is possible to replace the inequalities in
(\ref{*A}) with an asymptotic expansion. The critical values for $DS_n$ are 
$c_{jkn}(\alp)=\sqrt{2 \log \frac{en}{k-j}}+\kappa_n(\alp)$, while the above approximation
to the SAC penalty shows that for small intervals $k-j \ll n$ the critical values of $SAC_n$ 
are $\approx \sqrt{2 \log \frac{en}{k-j}} +\sqrt{2}\frac{\log \log (e(k-j))}{\sqrt{\log 
\frac{en}{k-j}}} + \tilde{\kappa}_n$ with $ \tilde{\kappa}_n <\kappa_n$. Plugging these
critical values in the expansion for the smallest detectable $\sqrt{|I_n|} \,\mu_{min}(n,|I_n|)$  shows 
that for $I_n \ll  n$ we obtain $e_n(|I_n|) \approx 1+\frac{b_{DS}}{\sqrt{ \log \frac{n}{|I_n|}}}$
$+o\Bl(\frac{1}{\sqrt{ \log \frac{n}{|I_n|}}}\Br)$ in the case of $DS_n$, and likewise
for $SAC_n$. Importantly, $b_{DS}>b_{SAC}$ and $b_{DS}$ is also larger than the corresponding 
constants for the blocked scan and the Bonferroni scan. Hence such an expansion would provide
a theoretical explanation of the superior performance of these three calibrations when
compared to $DS_n$ and therefore complement existing optimality results that are not sensitive 
enough to discern this effect, as was explained in Section~\ref{optimality}. We leave the
rigorous demonstration of such an optimality theory open.

Continuing with the proof we next check (\ref{*R}) for the blocked scan. By the union bound
\begin{align*}
\alp & = \Pr \left(\bigcup_{B=1}^{B_{max}} \Bl\{ \max_{(j,k]\in \mbox{Bth block}} T_{jk}(\Zn)
  > c_{B,n}\Bigl(\frac{\tilde{\alp}}{B}\Br)\Br\}\right)\\
& \leq \sum_{B=1}^{B_{max}} \frac{\tilde{\alp}}{B} \ \leq \ \tilde{\alp} (\log B_{max} +1)\
 \leq \ \tilde{\alp} (\log \log n +2),
\end{align*}
so 
\be \label{*C}
\frac{\tilde{\alp}}{B} \ \geq \ \frac{\alp}{(\log_2 n) (\log \log n +2)}.
\ee
In the case $B \geq 2$ the proofs of Theorems 6.1 and 2.1 in D\"{u}mbgen and Spokoiny~(2001)
show that for some constant $C$ (which is universal in this context) and for every $S\geq 1$:
\begin{align}
\Pr \Bl(\max_{(j,k]\in \mbox{Bth block}} T_{jk}(\Zn) & > \sqrt{2 \log \frac{en}{2^{B-1+s_n}} +
  S \log \log \frac{e^en}{2^{B-1+s_n}} } \Br) \nn \\
& \leq C \exp \Bl( (C-\frac{S}{C}) \log \log \frac{e^en}{2^{B-1+s_n}} \Br) \nn \\
& \leq C \Bl( \log \frac{e^en}{2^{B-1+s_n}} \Br)^{-2} \ \ \mbox{ by choosing $S \geq C(C+2)$}. 
\label{*D}
\end{align}
In the case $B=1$ the inequality (\ref{*D}) also holds for $S \geq 8$ as can be checked by
applying the union bound, the Gaussian tail bound and the fact that there are not more than
$n2^{s_n} \leq 2 n \log n$ intervals in the first block. Since we need to establish (\ref{*R})
only for interval lengths $k-j \leq n^p$, we need only consider block indices $B \leq \log_2 n^p
-s_n +2$. For those $B$, (\ref{*D}) is not larger than
$$
C \Bl( \log \frac{e^en}{2n^p} \Br)^{-2} \ \leq \ C \Bl( (1-p) \log n \Br)^{-2} \ \leq \ 
\frac{\alp}{(\log_2 n) (\log \log n +2)}\ \ \mbox{ eventually}.
$$
Comparing with (\ref{*C}) we immediately obtain
$$
c_{B,n}\Bigl(\frac{\tilde{\alp}}{B}\Br) \ \leq \ \sqrt{2 \log \frac{en}{2^{B-1+s_n}} +
S \log \log \frac{e^en}{2^{B-1+s_n}} }.
$$
Now (\ref{*R}) follows, because every interval $(j,k]$ belongs to a block with some index $B$.
Then $k-j < 2^{B-1+s_n}$ and so the critical value $c_{B,n}\Bigl(\frac{\tilde{\alp}}{B}\Br)$
for $T_{jk}(\Yn)$ is not larger than
$$
\sqrt{2 \log \frac{en}{k-j} +S \log \log \frac{e^en}{k-j}} \ =\ 
\sqrt{2 \log \frac{en}{k-j}} +O(1).
$$

Next we will establish the claim for the four calibrations that use the approximating set
of intervals. The sparsity of the approximating set makes it straightforward to establish
(\ref{*R}) for the null distribution, as will be seen below. On the other hand, we now have
to account for the approximation error incurred by not being able to match the support
of the signal exactly. The approximating set is constructed such that there is a bound on the
 error that is of the size needed:

\begin{Proposition} \label{approx}
For every $I = (j,k] \subset (0,n]$ there exists an interval $J$ in the approximating collection   
such that
$$
\frac{|I \cap J|}{\sqrt{|I||J|}}\ \geq \ \sqrt{1-\frac{2}{\sqrt{2\log \frac{en}{|I|\wedge |J|}}}}
\ \geq \ 1-\frac{1}{\sqrt{2\log \frac{en}{|I|\wedge |J|}}} -\frac{1}{\log \frac{en}{|I|\wedge |J|}}.
$$
\end{Proposition}

The proof of the proposition is below. 
So if $f_n(i)=\mu_n {\bs 1}(i \in I_n)$ is a signal from model (\ref{model}), then there exists 
an interval $J_n$ in the approximating set such that
$$
T_{J_n}(\Yn) \ = \ \frac{|I_n \cap J_n|}{\sqrt{|J_n|}} \mu_n + T_{J_n}(\Zn)\ 
=\ \frac{|I_n \cap J_n|}{\sqrt{|J_n||I_n|}} \sqrt{|I_n|} \,\mu_n +T_{J_n}(\Zn)
$$
by (\ref{a}).
So if $\sqrt{|I_n|} \,\mu_n =\sqrt{ 2\log \frac{en}{|I_n|}}+b_n$, 
then we get with Proposition~\ref{approx}
\begin{align}
T_{J_n}(\Yn) & \geq \Bl(1-\frac{1}{\sqrt{2\log \frac{en}{|I_n|}}}-\frac{1}{\log \frac{en}{|I_n|}}\Br)
  \Bl(\sqrt{ 2\log \frac{en}{|I_n|}}+b_n\Br) +T_{J_n}(\Zn)  \nn \\
& \geq \sqrt{ 2\log \frac{en}{|I_n|}} +b_n\Bl(1-\frac{1}{\sqrt{(1-p)\log (en)}}\Br)-1
-\sqrt{\frac{2}{(1-p)\log (en)}}+T_{J_n}(\Zn) 
  \ \ \ \mbox{ since $|I_n| \leq n^p$} \nn \\
& \geq \sqrt{ 2\log \frac{en}{|J_n|}} +b_n\Bl(1-\frac{1}{\sqrt{(1-p)\log (en)}}\Br)-1
  -\frac{\sqrt{2} \log (2e)}{(1-p) \sqrt{\log (en)}} +T_{J_n}(\Zn) \ \ 
  \mbox{ as $\frac{|I_n|}{|J_n|} \in [\frac{1}{2},2]$.} \label{compare2}
\end{align}
Now we can proceed as before: 
Suppose a generic test on the approximating set
assigns critical values $\tilde{c}_{jkn}(\alp)$ to the $T_{jk}(\Yn)$ that satisfy
\be  \label{*R*}
\tilde{c}_{jkn}(\alp) - \sqrt{2 \log \frac{en}{k-j}} \ \leq\ \tilde{b} \ \ 
\mbox{ for all $(j,k]$ in the approximating set
with $k-j \leq n^p$}.
\ee
Set $b_n$ such that the sum of the middle three terms in (\ref{compare2}) equals $\tilde{b}+z_{0.2}$.
Then, denoting the endpoints of $J_n$ by $J_n=(j_n,k_n]$:
$$
\Pr_{f_n}\Bl( \mbox{test rejects}\Br) \geq \Pr_{f_n}\Bl( T_{J_n}(\Yn)>\tilde{c}_{j_nk_nn}(\alp)\Br)
\geq \Pr\Bl(\sqrt{ 2\log \frac{en}{|J_n|}} +\tilde{b} +z_{0.2} +T_{J_n}(\Zn) >
  \tilde{c}_{j_nk_nn}(\alp)\Br) \geq 80\%
$$ 
by (\ref{*R*}). Hence the smallest detectable $\sqrt{|I_n|} \,\mu_{min}(n,|I_n|)$ for this generic test 
is not larger than $\sqrt{ 2\log \frac{en}{|I_n|}}+b_n$.
 Since the above choice for $b_n$ implies $b_n \ra \tilde{b} +1+ z_{0.2}$,
the claim of the theorem will follow as before upon verifying (\ref{*R*}). 

As an aside we note that in the case without an approximating set the bound on the smallest 
detectable $\sqrt{|I_n|} \,\mu_{min}(n,|I_n|)$ has the term $b+z_{0.2}$ in place of $b_n \approx
\tilde{b} +1 +z_{0.2}$, with $b$ from (\ref{*R}) and $\tilde{b}$ from (\ref{*R*}). 
Hence using the approximating set adds 1 to that bound but allows
to use $\tilde{b}$, and $\tilde{b} \leq b$ because the there are fewer intervals in the
approximating set to maximize over. We found in our simulations that using the approximating
set typically results in slightly more power.

(\ref{*R*}) clearly holds for $DS_n^{app}$, $SAC_n^{app}$ and {\sl BlockedScan$_n^{app}$}
since (\ref{*R}) holds for their counterparts $DS_n$, $SAC_n$ and {\sl BlockedScan$_n$},
and the former statistics cannot be larger than their counterparts since they maximize
over a subset of the intervals that their counterparts maximize over.
Therefore the claim of the theorem follows for these three calibrations. 

We note that we
established (\ref{*R}) by appealing to Theorem 2.1 in D\"{u}mbgen
and Spokoiny~(2001), which rests on their Theorem~6.1. The assumptions of that theorem
are difficult to check in general, as exemplified by the proof of their Theorem 2.1.
The sparsity of the approximating set makes it possible to alternatively establish (\ref{*R*})
directly with a simple application of the union bound and the Gaussian tail bound, similarly
to the following derivation for the Bonferroni scan:

$w_B := 2^{B-1+s_n}$ is an upper bound on the interval length in the $B$th block. It
follows from (\ref{count}) that
$$
\# (\mbox{$B$th block})\leq \left\{ \begin{array}{ll}
  \frac{n}{w_B} 8 \Bl( \log (en)\Br)^2 & \mbox{ if $B=1$}\\
  \frac{n}{w_B} 4 \log \Bl(\frac{n}{w_B} 2e\Br) & \mbox{ if $B\geq 2$}
  \end{array}
\right.
$$
Hence the Gaussian tail bound gives
\begin{align*}
v_n & \left(\frac{\alp}{\# (B\mbox{th block})\, B \,\sum_{i=1}^{B_{max}} \frac{1}{i}} \right)
 = \Phi^{-1}\Bl(\frac{\alp}{\# (B\mbox{th block})\, B \,\sum_{i=1}^{B_{max}} \frac{1}{i}}\Br)\\
& \leq \sqrt{2 \log \frac{\# (B\mbox{th block})\, B \,\sum_{i=1}^{B_{max}} \frac{1}{i}}{\alp}}\\
& \leq \sqrt{2 \log \frac{n}{w_B} + 2 \log \frac{8 (\log en)^2 B \sum_i \frac{1}{i}}{\alp}}  \\
& \leq \sqrt{2 \log \frac{n}{w_B} + 2 \log \frac{12 (\log en)^3 (\log \log n + \frac{3}{2})}{
  \alp}}  \\
& \ \ \mbox{ since } B \leq B_{max}\leq \frac{3}{2}\log n \mbox{ and } \sum_{i=1}^{B_{max}} 
  \frac{1}{i} \leq \log B_{max} + \frac{1}{2B_{max}}+0.58 \leq \log \log n + \frac{3}{2}  \\
& \leq \sqrt{2 \log \frac{n}{w_B}} + \frac{\log \Bl(\frac{12}{\alp} (\log n)^3 (\log \log n
  + \frac{3}{2})\Br)}{\sqrt{2 \log \frac{n}{w_B}}}  \\
& \leq \sqrt{2 \log \frac{n}{w_B}} +O(1)\ \ \ \mbox{ if $w_B \leq n^p$,\ $p<1$}
\end{align*}
which establishes (\ref{*R*}) for the Bonferroni scan. 

It remains to prove Proposition~\ref{approx}: Let $\ell$ be such that $|I| \in [m_{\ell},
2m_{\ell})$. Elementary considerations show that one can pick $J \in {\cal J}_{\ell}$ such
that $|I \triangle J| \leq d_{\ell}$. Hence
\begin{align*}
\frac{|I \cap J|}{\sqrt{|I||J|}} & = \sqrt{\frac{|I|-|I\setminus J|}{|I|}} 
   \sqrt{\frac{|J|-|J\setminus I|}{|J|}}\\
& = \sqrt{1-\alp \frac{|I\triangle J|}{|I|}} \sqrt{1-(1-\alp) \frac{|I\triangle J|}{|J|}}
  \ \ \ \mbox{ where } \alp=\frac{|I \setminus J|}{|I\triangle J|}\\
& \geq \sqrt{1-\frac{|I\triangle J|}{\min (|I|,|J|)}}\ \ \mbox{ since }
  (1-\alp x)(1-(1-\alp)x) \geq 1-x\\
& \geq \sqrt{1-\frac{d_{\ell}}{m_{\ell}}}
\end{align*}
If $m_{\ell} > \sqrt{2 \log \frac{en}{m_{\ell}}}$, then 
$$
\frac{d_{\ell}}{m_{\ell}}
\leq \frac{\frac{m_{\ell}}{\sqrt{2 \log \frac{en}{m_{\ell}}}} +1}{m_{\ell}}
\leq \frac{2}{\sqrt{2 \log \frac{en}{|I| \wedge |J|}}}
$$
 and the first inequality of the claim follows. 
If $m_{\ell} \leq \sqrt{2 \log \frac{en}{m_{\ell}}}$, then
$d_{\ell}= \Big\lceil \frac{m_{\ell}}{\sqrt{2 \log \frac{en}{m_{\ell}}}} \Big\rceil =1$, so
by the definition of ${\cal J}_{\ell}$ we can take $J:=I$ and the first inequality of
the claim also holds. The second inequality of the proposition follows from 
$\sqrt{1-x} \geq 1-\frac{1}{2}x-\frac{1}{2}x^2$ for $x \in (0,1)$. $\Box$

The proof of the theorem uses the following bounds on the number of intervals in the approximating
set: There are no more than $\frac{n}{d_{\ell}}$ possible left endpoints for intervals
in ${\cal J}_{\ell}$ and for each left endpoint there are no more than 
$\frac{m_{\ell}}{d_{\ell}}$ right endpoints, hence
$$
\# {\cal J}_{\ell} \ \leq \ \frac{n m_{\ell}}{d_{\ell} d_{\ell}} \ \leq \ 
\frac{n m_{\ell} 2 \log \frac{en}{m_{\ell}}}{m_{\ell} m_{\ell}}\ =\ 2n 2^{-\ell} 
\log(e2^{-\ell}n).
$$
Therefore
\be \label{count}
\# (\mbox{$B$th block}) \left\{ \begin{array}{ll}
  \leq \ \sum_{\ell=0}^{s_n-1} 2n 2^{-\ell} \log(e2^{-\ell}n)\ \leq \ 4n 
    \log (en) & \mbox{ if $B=1$} \\
  \ =\ \# {\cal J}_{B-2+s_n}\ \leq \ 2n 2^{-B+2-s_n} \log (e2^{-B+2-s_n}n)
    \ \leq \ 8n 2^{-B}  & \mbox{ if $B\geq 2$}  
  \end{array}
\right.
\ee
since $2^{-s_n} \leq \frac{1}{\log n}$. $\Box$

\subsection*{References}

\begin{description}
\item[] Arias-Castro, E., Donoho, D. and Huo, X. (2005). Near-optimal detection of geometric
objects by fast multiscale methods. {\sl IEEE Trans. Inform. Theory} {\bf 51}, 2402--2425.
\item[] Arias-Castro, E., Cand\`{e}s, E. and Durand, A. (2011). Detection of an
anomalous cluster in a network. {\sl Ann. Statist.} {\bf 39}, 278--304.
\item[] Arias-Castro, E., Castro, R.M., T\'{a}nczos, E. and Wang, M. (2018). Distribution-free
detection of structured anomalies: Permutation and rank-based scans.
{\sl J. Amer. Statist. Assoc.} {\bf113}, 789--801.
\item Bardenet, R. and Maillard, O. (1995). Concentration inequalities for sampling
without replacement. {\sl Bernoulli} {\bf 21}, 1361-1385.
\item Barndorff-Nielsen, O.E. (1986). Inference on full or partial parameters based on the
standardized signed log likelihood ratio. {\sl Biometrika} {\bf 73}, 307--322.
\item[] Bentkus, V. K. and Dzindzalieta, D. (2015). A tight Gaussian bound for weighted sums 
of Rademacher random variables. {\sl Bernoulli} {\bf 21}, 1231–1237.
\item[] Broemeling, L.D. and Tsurumi, H. (1987). {\sl Econometrics and Structural Change.}
Marcel Dekker, New York.
\item[] Brown, L.B. and Low, M.G. (1996). Asymptotic equivalence of nonparametric regression
and white noise. {\sl Ann. Statist.} {\bf 24}, 2384--2398.
\item[] Chan, H.P. and Lai, T.L. (2006). Maxima of asymptotically gaussian random fields and
moderate deviation approximations to boundary crossing probabilities of sums of random
variables with multidimensional indices. {\sl Ann. Probab.} {\bf 34}, 80--121.
\item[] Chan, H.P. and Walther, G. (2013). Detection with the scan and the average likelihood
ratio. {\sl Statistica Sinica} {\bf 23}, 409--428.
\item[] Datta, P. and Sen, B. (2018). Optimal inference with a multidimensional multiscale statistic.
arXiv:1806.02194
\item[] Darling, D.A. and Erd\"{o}s, P. (1956). A limit theorem for the maximum of normalized
sums of independent random variables. {\sl Duke Math. J.} {\bf 23}, 143--154.
\item[] de la Pe\~{n}a, V.H., Lai, T.L. and Shao, Q.M. (2009). {\sl Self-Normalized Processes:
Theory and Statistical Applications.} Springer, Berlin.
\item[] D\"{u}mbgen, L. and Spokoiny, V.G. (2001). Multiscale testing of qualitative
hypotheses. {\sl Ann. Statist.} {\bf 29}, 124--152.
\item[] D\"{u}mbgen, L. (2002). Application of local rank tests to nonparametric regression.
{\sl J. Nonparametric Statistics} {\bf 14}, 511-537.
\item[] D\"{u}mbgen, L. and Walther, G.~(2008). Multiscale inference about a density.
{\sl Ann. Statist.} {\bf 36}, 1758--1785.
\item[] Enikeeva, F., Munk, A. and Werner, F. (2018). Bump detection in heterogeneous
Gaussian regression. {\sl Bernoulli} {\bf 24}, 1266--1306.
\item[] Frick, K., Munk, A. and Sieling, H. (2014). Multiscale change point inference.
{\sl J. R. Stat. Soc. Ser. B.} {\bf 76}, 495--580.
\item[] Glaz, J., Naus, J. and Wallenstein, S. {\sl Scan Statistics.} Springer Series in Statistics.
Springer-Verlag, New York, 2001. 
\item[] Hoeffding, W. (1952). The large-sample power of tests based on permutations of 
observations. {\sl Ann. Math. Statist.} {\bf 23}, 169--192.
\item[] Hogan, M.L. and Siegmund D. (1986). Large deviations for the maxima of some random
fields. {\sl Adv. Appl. Math.} {\bf 7}, 2--22.
\item[] Jung, I. and Cho, H. (2015). A nonparametric spatial scan statistic for continuous data.
{\sl International Journal of Health Geographics} {\bf 14}, 30.
\item[] Kabluchko, Z. (2011). Extremes of the standardized gaussian noise. {\sl Stochastic Processes 
and their Applications} {\bf  121}, 515--533. 
\item[] K\"{o}nig, C., Munk, A. and Werner, F. (2020). Multidimensional
multiscale scanning in exponential families: Limit theory and statistical consequences.
{\sl Ann. Statist.} {\bf 48}, 655-678.
\item[] Kulldorff, M. (1997). A spatial scan statistic. {\sl Comm. Statist. Theory Methods}
{\bf 26}, 1481–1496.
\item[] Levin, B. and Kline, J. (1985). The CUSUM test of homogeneity with an application in
spontaneous abortion epidemiology. {\sl Statistics in Medicine} {\bf 4}, 469--488.
\item[] Loader, C. R. (1991). Large-deviation approximations to the distribution of the scan statistics.
{\sl Adv. Appl. Prob.} {\bf 23}, 751--771.
\item[] Naus, J. I. and Wallenstein, S. (2004). Multiple window and cluster size scan procedures.
{\sl Meth. Comp. Appl. Probab.} {\bf 6}, 389--400.
\item[] Neill, D. and Moore, A. (2004a). A fast multi-resolution method for detection of significant spatial 
disease clusters. {\sl Adv. Neural Inf. Process. Syst.} {\bf 10}, 651–658.
\item[] Neill, D. and Moore, A. (2004b). Rapid detection of significant spatial disease clusters. In {\sl Proc. 
Tenth ACM SIGKDD International Conference on Knowledge Discovery and Data Mining} 256– 265. ACM, New York.
\item[] Olshen, A.B., Venkatraman, E.S., Lucito, R. and Wigler, M. (2004). Circular binary segmentation for the 
analysis of array-based DNA copy number data. {\sl Biostatistics} {\bf 5}, 557–572.
\item[] Pinelis, I. (2012). An asymptotically Gaussian bound on the Rademacher tails.
{\sl Electron. J. Probab.} {\bf 17}, 1--22.
\item[] Proksch, K, Werner, F. and Munk, A. (2018). Multiscale scanning in inverse problems.
{\sl Ann. Statist.} {\bf 46}, 3569--3602.
\item[] Rivera,  C. and Walther, G. (2013). Optimal detection of a jump in the intensity of a Poisson
process or in a density with likelihood ratio statistics. {\sl Scand. J. Stat.} {\bf 40}, 752-769.
\item[] Rohde, A. (2008). Adaptive goodness-of-fit tests based on signed ranks.
{\sl Ann. Statist.} {\bf 36}, 1346--1374.
\item[] Romano, J.P. (1989). Bootstrap and randomization tests of some nonparametric hypotheses.
{\sl Ann. Statist.} {\bf 17}, 141--159.
\item[] Sharpnack, J, and Arias-Castro, E. (2016). Exact asymptotics for the scan statistic and fast
alternatives. {\sl Elect. J. Statist.} {\bf 10}, 2641-2684.
\item[] Siegmund, D. (1986). Boundary crossing probabilities and statistical applications.
{\sl Ann. Statist.} {\bf 14}, 361--404.
\item[] Siegmund, D. (2017). Personal communication.
\item[] Siegmund, D. and Venkatraman, E.S. (1995). Using the generalized
likelihood ratio statistic for sequential detection of a change-point. {\sl Ann.
Statist.} {\bf 23}, 255--271.
\item[] Walther, G. (2010). Optimal and fast detection of spatial clusters with scan statistics.
{\sl Ann. Statist.} {\bf 38}, 1010--1033.
\item[] Walther, G. (2021). Tail bounds for empirically standardized sums. Manuscript.
\item[] Yao, Q. (1993a). Tests for change-points with epidemic alternatives.
{\sl Biometrika} {\bf 80}, 179--191.
\item[] Yao, Q. (1993b). Boundary-crossing probabilities of some random fields related to 
likelihood ratio tests for epidemic alternatives. {J. Appl. Prob.} {\bf 30}, 52--65.
\item[] Yu, Y. (2020). A review of minimax rates in change point detection and localisation.
arXiv preprint arXiv:2011.01857.
\end{description}

\end{document}